\documentclass[pdflatex,sn-mathphys-num,12pt]{sn-jnl}

\usepackage{amsmath,amsfonts,amssymb}
\usepackage{amsthm}
\usepackage{graphicx}
\usepackage{subcaption}
\usepackage{booktabs}
\usepackage{enumitem}
\usepackage{algorithm}
\usepackage{algorithmicx}
\usepackage[noend]{algpseudocode}
\algrenewcommand\algorithmicrequire{\textbf{Input:}}
\algrenewcommand\algorithmicensure{\textbf{Output:}}
\usepackage[title]{appendix}
\usepackage{xcolor}
\usepackage{placeins}
\theoremstyle{thmstyleone}
\newtheorem{teo}{Theorem}
\newtheorem{lem}[teo]{Lemma}
\newtheorem{prop}[teo]{Proposition}
\newtheorem{cor}[teo]{Corollary}

\theoremstyle{thmstyletwo}
\newtheorem{rem}{Remark}
\newtheorem{example}{Example}

\theoremstyle{thmstylethree}
\newtheorem{assumption}{Assumption}

\newcommand{\R}{\mathbb{R}}
\newcommand{\scalar}[2]{\langle #1 , #2 \rangle}
\DeclareMathOperator*{\argmin}{arg\,min}
\DeclareMathOperator{\prox}{prox}
\DeclareMathOperator*{\proj}{proj}

\setlist[enumerate]{leftmargin=.5in}
\setlist[itemize]{leftmargin=.5in}

\begin{document}

\title[DC optimization for quasinormed TV image reconstruction]
      {DC optimization for quasinormed total variation image reconstruction}

\author[1,2]{\fnm{Pedro} \sur{Merino}}\email{pedro.merino@epn.edu.ec}

\author[1,2]{\fnm{Daniela} \sur{Riera}}\email{daniela.riera@epn.edu.ec}

\author*[3]{\fnm{David} \sur{Villac\'{\i}s}}\email{davillacis@uloyola.es}

\affil[1]{\orgname{Research Center for Mathematical Modeling (MODEMAT)},
  \orgaddress{\city{Quito}, \country{Ecuador}}}

\affil[2]{\orgdiv{Departamento de Matem\'atica},
  \orgname{Escuela Polit\'ecnica Nacional},
  \orgaddress{\city{Quito}, \country{Ecuador}}}

\affil*[3]{\orgdiv{Departamento de M\'etodos Cuantitativos},
  \orgname{Universidad Loyola Andaluc\'ia},
  \orgaddress{\city{Seville}, \country{Spain}}}

\abstract{We address image reconstruction regularized by the quasinormed total variation $\text{TV}^q$, $q\in(0,1)$, with a general linear forward operator, and its numerical solution by difference-of-convex (DC) programming. Since $t^q$ has unbounded slope at the origin, the objective admits no DC splitting whose convex part is a total variation functional. For concave penalties of finite slope, which include a Huber-type smoothing of $\text{TV}^q$, we construct a splitting whose convex part is a total variation functional with a single scalar weight and whose second part is continuously differentiable, and we bound the smoothing error. The subproblems are solved by a primal--dual method that returns a computable duality gap. Stopping on this gap certifies the descent of every outer step, and every limit point is a critical point of the smoothed problem. Experiments on denoising, deblurring, inpainting and tomography compare the method with iteratively reweighted $\ell_1$ and a semismooth Newton method.}



\maketitle

\section{Introduction}
\label{sec:introduction}
Image reconstruction problems are often ill-posed inverse problems, and their variational
treatment requires a balance between a data fidelity term and a regularization term. Total
variation is the standard choice of regularizer when edges have to be preserved. The nonconvex variant $\text{TV}^q$ with $q\in(0,1)$  has been used to enhance edge preservation  by promoting sparsity on the image gradient. In this work, we are interested in the numerical solution of the problem:
\begin{equation}\label{eq:tvq_problem}
  \min_{u}\;\frac12\|Au-f\|^2+\beta\sum_{i=1}^m\|(\mathbb{K}u)_i\|^q ,
\end{equation}
where $u\in\mathbb{R}^n$ is the vectorized image, $f\in\mathbb{R}^M$ the observed data, $A:\mathbb{R}^n\to\mathbb{R}^M$ a bounded linear forward operator, e.g., the identity for denoising, a convolution for deblurring, a coordinate projection for inpainting, or an X-ray transform for tomography, and $\mathbb{K}$ the discrete gradient, so that $(\mathbb{K}u)_i\in\mathbb{R}^2$ is the gradient at pixel $i$ and $\|(\mathbb{K}u)_i\|$ its magnitude. The parameter $\beta>0$ weights the regularization, and the penalty $\sum_{i=1}^m\|(\mathbb{K}u)_i\|^q$ in \eqref{eq:tvq_problem} is referred to as the quasinormed total variation regularization.

For $0<q<1$, the regularizer is non-convex, non-differentiable with unbounded slope at the origin. The last property is the reason why many algorithms in the literature that linearize the penalty work with a surrogate instead of $\sum_i\|(\mathbb{K}u)_i\|^q$ itself. The iteratively reweighted $\ell_1$ scheme (IRL1) of Ochs et al.\ \cite{ochs2015iteratively} is applied to the shifted power $(t+\varsigma)^q$, $\varsigma>0$, whose slope at the origin is $q\varsigma^{q-1}$. The semismooth Newton method of Hintermüller and Wu \cite{hintermuller2013nonconvex} replaces $t^q$ near the origin by a quadratic function, so that the curvature of the penalty stays bounded. In both cases, the function actually minimized has \emph{finite} slope at the origin.

Nonconvex gradient penalties go back to the half-quadratic models of Geman and Reynolds \cite{geman1992constrained}. Power penalties with $q<1$ were proposed for sparse recovery by Chartrand \cite{chartrand2007exact}, and Nikolova et al.\ \cite{nikolova2010fast} developed fast methods for nonsmooth nonconvex regularization in image restoration. Hintermüller et al.\ \cite{hin_valk} showed that in function space the $\text{TV}^q$ energy degenerates, which is one reason to pose the problem in the discrete setting, as we do. Lanza et al.\ \cite{lanza2016convex} choose the parameter of a concave penalty so that the whole objective remains convex, and the minimax concave penalty has been applied to image restoration in \cite{you2019mcp}. The same penalties appear in statistics as SCAD \cite{fan2001variable} and MCP \cite{zhang2010nearly}, and reweighted $\ell_1$ minimization \cite{candes2008enhancing} is the standard way of linearizing them. A difference-of-convex structure underlies all of these schemes. DC programming and the DC algorithm are surveyed in \cite{taole,hor,lethi2018dc}; boosted variants add a line search along the DCA step \cite{aragon2018accelerating,aragon2020boosted} and have been applied to image restoration under Rician noise with a smoothed total variation term \cite{wu2022rician}. Smoothing of nonsmooth nonconvex penalties as a general tool is studied by Chen \cite{chen2012smoothing}. The splitting proposed in this paper differs from the above in that the convex part is itself a nonsmooth total variation problem, which we solve inexactly to a certified duality gap, and in that the smoothing is local, so that the DC remainder vanishes away from the edges.

The methods of Ochs et al.\ and of Hintermüller and Wu, and the splitting method of Fang et al.\ \cite{fang2022nonconvex}, handle the non-convexity in different ways, each with its own limitations. IRL1 majorizes the concave penalty by its tangent at the current iterate, so that every step solves a total variation problem with a spatially varying weight, which changes at every iteration. When this subproblem is solved through a linear system, as in iteratively reweighted least squares or Newton-type methods, the weights enter the coefficient matrix of that linear system, and its conditioning degrades as they spread across the image. In a primal--dual solver, by contrast, the radii of the dual constraints are the only change. Hintermüller and Wu solve their smoothed problem by a globalized semismooth Newton method with a trust-region regularization; a regularization of the generalized Hessian, with weight in $[1-q,\,2-q]$, keeps the modified Newton matrix positive definite \cite[Theorem~3.2]{hintermuller2013nonconvex}. Fang et al.\ consider $\text{TV}^q$ with an $\ell_1$ fidelity term and apply a majorization--minimization scheme within an ADMM framework, which decouples the regularization through the split variable $v=\mathbb{K}u$ and absorbs the non-convexity into a separable blockwise proximal operator. This avoids a spatially varying weight in the subproblem, but their convergence result requires the ADMM penalty factor $\varrho$ to be sufficiently large \cite[Theorem~1]{fang2022nonconvex}, while large values of $\varrho$ make the method slow, so that $\varrho$ is increased along the iteration.

\subsection{Contributions}
\label{sec:contributions}
In this article, we propose a novel difference-of-convex algorithm (DCA) for $\text{TV}^q$. We replace $t^q$ by a linear function on a small interval at the origin and keep it, up to a shift, elsewhere. The penalty is then the sum of a total variation term with a scalar weight and a concave function that vanishes on the flat regions of the image. We keep the first as the convex part of the splitting and linearize the second, so that every subproblem is the same total variation problem perturbed by a linear term that acts only at the edges. In IRL1, by contrast, the whole penalty is linearized and the subproblem is a total variation problem with a weight that varies over the image. The smoothing comes from \cite{merino}, where it was used for sparse optimal control problems with a nonconvex cost. The present setting, a penalty on the image gradient under a general forward operator, is what leads to the admissible class, to the result that no such splitting exists without smoothing, and to the primal--dual solution of the subproblems with a computable certificate. Our main contributions are summarised below. 

\begin{itemize}
\item \emph{A DCA with a fixed total variation subproblem.} In our algorithm, every subproblem is the same total variation problem with one scalar weight, perturbed by a linear term that acts only at the edges, unlike IRL1, whose weight varies over the image and at every step. This requires a surrogate of $t^q$ with finite slope, which is unavoidable (Corollary~\ref{cor:no-splitting}), and the concave part is then $C^1$ with a closed-form gradient (Lemma~\ref{lem:H_dif}). 
The subproblems are solved by a primal--dual method in which $A$ is only applied and never inverted, so the algorithm applies to any linear forward operator.
\item \emph{A certified algorithm.} The primal--dual method returns a computable duality gap (Proposition~\ref{prop:gap}) that bounds a stationarity measure and sets the inner stopping rule, so no sequence of tolerances is chosen in advance. Each outer step decreases the energy by a fixed fraction of that measure. For a positive tolerance the number of outer iterations is bounded, and for zero tolerance every limit point is a critical point of the smoothed problem (Proposition~\ref{prop:descent}, Theorem~\ref{teo:convergence}).
\item \emph{A competitive and efficient method.} On denoising, deblurring, inpainting and tomography, the algorithm brings the energy within $1\%$ of the best value found on $116$ of $130$ runs, against $107$ for IRL1, and matches the reconstruction quality of Hinterm\"uller and Wu \cite{hintermuller2013nonconvex}, with less work on nine of the thirteen problems (Section~\ref{sec:experiments}).
\end{itemize}

The paper is organized as follows. Section~\ref{sec:problem} states the problem and the class of admissible penalties. Section~\ref{sec:dcsplit} constructs the splitting, and Section~\ref{sec:smoothing_error} bounds the error of the smoothing. Section~\ref{sec:pdhg} derives the algorithm, its duality gap and its convergence. Section~\ref{sec:experiments} contains the numerical study, and Section~\ref{sec:conclusions} the conclusions and limitations. The appendices collect the proofs, the implementation of the comparison methods and reproducibility details.

\bmhead{Use of AI tools} The conception of the method, including the DC splitting, its analysis, the algorithm and the design of the experiments, is the work of the authors. The authors implemented the core of the algorithm; starting from this implementation, Claude Opus 5.5 and Claude Fable 5.2 (Anthropic) were used to extend the code of the methods and the experiments, and to improve the presentation of the proofs and the text (Appendix~\ref{app:repro}). The authors verified the code and the proofs and take full responsibility for the content of this publication.

\subsection{Notation and preliminaries}
Throughout, $u\in\mathbb{R}^n$ is the vectorized image, $f\in\mathbb{R}^M$ the observed data and $A:\mathbb{R}^n\to\mathbb{R}^M$ the linear forward operator; $\|\cdot\|$ is the Euclidean norm on the corresponding space. The discrete gradient $\mathbb{K}:\mathbb{R}^n\to\mathbb{R}^{2m}$ stacks the horizontal and vertical forward differences, $\mathbb{K}=[\mathbb{K}_x;\mathbb{K}_y]$, set to zero across the last column and the last row, so that $m=n$ and $(\mathbb{K}u)_i=(\mathbb{K}_x^iu,\mathbb{K}_y^iu)\in\mathbb{R}^2$ is the gradient at pixel $i$, where $\mathbb{K}_x^i$ and $\mathbb{K}_y^i$ are the $i$-th rows of $\mathbb{K}_x$ and $\mathbb{K}_y$. We write $\phi_i(u):=\|(\mathbb{K}u)_i\|$ for its magnitude, $i=1,\dots,m$. 

A function $J$ is a difference of convex (DC) function if $J=G-H$ with $G$ proper, convex and lower semicontinuous and $H$ convex and finite; $(G,H)$ is a DC splitting of $J$. A point $u$ is a critical point of $G-H$ if $\partial G(u)\cap\partial H(u)\neq\emptyset$.

\section{Regularized TV$^q$ image reconstruction problem}
\label{sec:problem}
Problem \eqref{eq:tvq_problem} is not solved directly. Since $t\mapsto t^q$ has unbounded derivative at $t=0$, no linear function $\lambda t$ majorizes it, and we prove in Corollary~\ref{cor:no-splitting} that \eqref{eq:tvq_problem} admits no difference-of-convex splitting of the form used in this paper. We therefore replace $t^q$ by a surrogate. This section states the properties that the surrogate must satisfy and shows that the Huber-type smoothing introduced in \cite{merino}, which is the one used in our algorithm and experiments, has them (Lemma~\ref{lem:huber-admissible}). Given a function $\varphi:[0,\infty)\to\mathbb{R}$, which we call the penalty, we consider the restoration problem
{\renewcommand{\theHequation}{eqP}%
\begin{equation}
  \min_{u} J_\varphi(u):=\frac{1}{2}\|Au-f\|^2 + \beta\sum_{i=1}^m \varphi(\phi_i(u)).
  \tag{$P_{\varphi}$}
  \label{eq:Pphi}
\end{equation}}
The penalties we consider satisfy the following assumption.
\begin{assumption}
\label{ass:penalty}
A function $\varphi:[0,\infty)\to\mathbb{R}$ is an \emph{admissible penalty} if
\begin{enumerate}[label=(A\arabic*),leftmargin=.7in]
  \item $\varphi(0)=0$, $\varphi$ is continuous on $[0,\infty)$ and nondecreasing;
        \label{ass:A1}
  \item $\varphi$ is concave on $[0,\infty)$ and continuously differentiable on
        $(0,\infty)$;\label{ass:A2}
  \item the slope at the origin is finite,
        $\displaystyle\delta_\varphi:=\varphi'(0^+)=\lim_{t\downarrow0}\frac{\varphi(t)}{t}<\infty$.
        \label{ass:A3}
\end{enumerate}
\end{assumption}

\begin{example}[Total variation] 
The identity $\varphi(t)=t$ is admissible, with $\delta_\varphi=1$, and \eqref{eq:Pphi} is then the total variation problem, the case $q=1$ of \eqref{eq:tvq_problem}. It is the baseline of the experiments in Section~\ref{sec:experiments}. We distinguish the total variation with a \emph{scalar weight}, $\lambda\sum_{i=1}^m\phi_i(u)$ with $\lambda>0$, which penalizes every pixel equally, from the \emph{spatially weighted} total variation $\sum_{i=1}^m\alpha_i\phi_i(u)$ with a vector of weights $\alpha\in\mathbb{R}^m_{+}$, which is the subproblem of iteratively reweighted schemes. The subproblems of our algorithm are of the first kind.
\end{example}

As a consequence of Assumption~\ref{ass:penalty}, an admissible penalty $\varphi$ is nondecreasing and concave, with finite slope at the origin. The next lemma collects some further properties that we will use in the sequel.

\begin{lem}\label{lem:varphi_prime_nneg_ninc}
Let $\varphi$ be an admissible penalty. Then:
\begin{enumerate}[label=(\roman*),leftmargin=.4in]
  \item $\varphi'$ is nonnegative and nonincreasing on $(0,\infty)$;
  \item $\lim_{t\downarrow0}\varphi'(t)=\sup_{t>0}\varphi'(t)=\delta_\varphi$; in particular,
        $\delta_\varphi\ge0$ and $\varphi'\le\delta_\varphi$ on $(0,\infty)$;
  \item $\varphi(t)\le\delta_\varphi t$ for every $t\ge0$.
\end{enumerate}
\end{lem}

\begin{proof}
\emph{(i)} It is clear that by assumptions \ref{ass:A1} and \ref{ass:A2} that $\varphi'$ is nonnegative. Furthermore, since $\varphi$ is concave by \ref{ass:A2}, its tangent lines at $0<s<t$ lie above its graph, which gives $\varphi'(t)\le\bigl(\varphi(t)-\varphi(s)\bigr)/(t-s)\le\varphi'(s)$.

\emph{(ii)} By (i), the limit $L:=\lim_{s\downarrow0}\varphi'(s)=\sup_{s>0}\varphi'(s)$ belongs to $[0,\infty]$. For $t>0$, $\varphi$ is continuous on $[0,t]$ and differentiable on $(0,t)$ and, since $\varphi(0)=0$, the mean value theorem gives $\xi_t\in(0,t)$ with $\varphi(t)/t=\varphi'(\xi_t)$. Since $\varphi'$ is nonincreasing and $\xi_t<t$,
\begin{equation}
  \varphi'(t)\;\le\;\frac{\varphi(t)}{t}\;\le\;L,\qquad t>0 .
  \label{eq:phi-slope}
\end{equation}
Letting $t\downarrow0$, \ref{ass:A3} gives $L\le\delta_\varphi\le L$, i.e. $L=\delta_\varphi$. The remaining claims follow from (i).

\emph{(iii)} For $t>0$, \eqref{eq:phi-slope} and (ii) give $\varphi(t)/t\le\delta_\varphi$. For $t=0$, both sides vanish by \ref{ass:A1}.
\end{proof}

\begin{figure}[ht]
  \centering
  \includegraphics[width=1.0\textwidth]{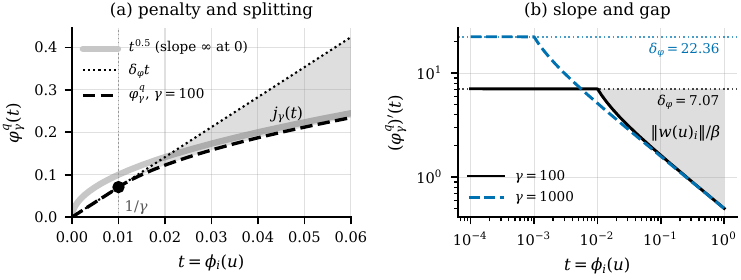}
  \caption{Huber-type smoothing \eqref{eq:reg} for $q=0.5$. \textbf{(a)} For $\gamma=100$, $\varphi_\gamma^q$ coincides with the line $\delta_\varphi t$ up to the knot $1/\gamma$, and the shaded gap between them is $j_\gamma$ of \eqref{j}, the function the algorithm linearizes. \textbf{(b)} The derivative for $\gamma=100$ and $1000$ is bounded by $\delta_\varphi=7.07$ and $22.36$, and the shaded gap $\delta_\varphi-(\varphi_\gamma^q)'$ for $\gamma=100$ is $\|w(u)_i\|/\beta$, see \eqref{eq:tilt-bound}.}
  \label{fig:huber_smoothing}
\end{figure}

The admissible penalty used in the algorithm of Section~\ref{sec:pdhg} and in the numerical experiments of Section~\ref{sec:experiments} utilizes the Huber-type local smoothing of $t^q$ introduced in \cite{merino}. It replaces $t^q$ by a linear function below a knot located at $1/\gamma$, and by a shifted version of $t^q$ above it, so that its slope at the origin is finite, see Figure~\ref{fig:huber_smoothing}. For $\gamma \gg 1$, we set
\begin{equation}
  \varphi_\gamma^q(t)=\left\{\begin{array}{ll}
    q^q\gamma^{1-q}\,t, & 0\le t\le 1\slash\gamma,\\[8pt]
    \Bigl(t+\dfrac{q-1}{\gamma}\Bigr)^{q}, & t>1\slash\gamma.
  \end{array}\right.
  \label{eq:reg}
\end{equation}

It is straightforward to see that $\varphi^q_\gamma \to |\cdot|^q$ as $\gamma \to \infty$. The next lemma establishes that $\varphi^q_\gamma$ is an admissible penalty.
\begin{lem}
\label{lem:huber-admissible} 
For $q\in(0,1)$ and $\gamma \gg 1$, the function $\varphi_\gamma^q$ of \eqref{eq:reg} is an admissible penalty, with slope at the origin $\delta_{\varphi}=q^q\gamma^{1-q}$. Moreover, $\varphi_\gamma^q$ is continuously differentiable on $[0,\infty)$, with
\begin{equation}
  (\varphi_\gamma^q)'(t)=\left\{\begin{array}{ll}
    \delta_\varphi, & 0\le t\le 1\slash\gamma,\\[6pt]
    q\Bigl(t+\dfrac{q-1}{\gamma}\Bigr)^{q-1}, & t>1\slash\gamma,
  \end{array}\right.
  \label{eq:reg-deriv}
\end{equation}
and $0<(\varphi_\gamma^q)'\le\delta_\varphi$.
\end{lem}

\begin{proof}
At $t=1/\gamma$, both branches of \eqref{eq:reg} take the value $q^q\gamma^{-q}=(q/\gamma)^q$, and both have slope $\delta_\varphi$, since $q\bigl(\frac1\gamma+\frac{q-1}{\gamma}\bigr)^{q-1}=q\bigl(\frac{q}{\gamma}\bigr)^{q-1} =q^q\gamma^{1-q}$. Hence $\varphi_\gamma^q$ is continuously differentiable on $[0,\infty)$, with derivative given in \eqref{eq:reg-deriv}. For $t>1/\gamma$, we have $t+\frac{q-1}{\gamma}>\frac{q}{\gamma}$ and $q-1<0$, so the second branch of \eqref{eq:reg-deriv} is positive, nonincreasing in $t$ and bounded by $q(q/\gamma)^{q-1}=\delta_\varphi$. Therefore, $(\varphi_\gamma^q)'$ is positive and nonincreasing on $[0,\infty)$, i.e., $\varphi_\gamma^q$ is nondecreasing and concave, and, since $\varphi_\gamma^q(0)=0$, \ref{ass:A1} and \ref{ass:A2} hold. Finally, $\varphi_\gamma^q(t)/t=\delta_\varphi$ for $0<t\le1/\gamma$, which gives \ref{ass:A3}.
\end{proof}

For this penalty, \eqref{eq:Pphi} is the regularized problem
{\renewcommand{\theHequation}{eqPgamma}%
\begin{equation}
  \min_{u \in \mathbb{R}^n} J_\gamma(u):= \frac{1}{2}\|Au-f\|^2
  + \beta \sum_{i=1}^m \varphi_{\gamma}^q\bigl(\phi_i(u)\bigr).
  \tag{\text{$P_{\gamma}$}}
  \label{eq:Pgamma}
\end{equation}}

\section{The DC splitting}
\label{sec:dcsplit}
A difference-of-convex decomposition of the objective $J_\varphi$ of \eqref{eq:Pphi} is the basis of our method. In this section we construct one for every admissible penalty, acting pixel by pixel on $\varphi$. Its first function is a total variation restoration functional with the scalar weight $\beta\delta_\varphi$, and its second function is continuously differentiable. For the penalty $t^q$ itself no such decomposition exists, which is why the smoothing of Section~\ref{sec:problem} is needed.

\begin{prop}
\label{prop:sufficiency}
Let $\varphi$ be admissible and put
\begin{equation}
  j_\varphi(t):=\delta_\varphi\,t-\varphi(t),
  \qquad t\ge0 ,
  \label{eq:j}
\end{equation}
so that $\varphi(t)=\delta_\varphi t-j_\varphi(t)$, and
\begin{equation}
  G_\varphi(u)=\frac12\|Au-f\|^2+\beta\delta_\varphi\sum_{i=1}^m \phi_i(u),
  \qquad
  H_\varphi(u)=\beta\sum_{i=1}^m j_\varphi\bigl(\phi_i(u)\bigr).
  \label{eq:GH}
\end{equation}
Then $G_\varphi$ and $H_\varphi$ are convex on $\mathbb{R}^n$, for every bounded linear
$A$ and every $\beta>0$, and
\begin{equation}
  J_\varphi=G_\varphi-H_\varphi .
  \label{eq:dc}
\end{equation}
We call $G_\varphi$ and $H_\varphi$ the DC components of $J_\varphi$ and $\lambda:=\beta\delta_\varphi$ the weight of the splitting. The component $G_\varphi$ is a total variation restoration functional that depends on $\varphi$ only through the scalar $\delta_\varphi$.
\end{prop}

\begin{proof}
The identity \eqref{eq:dc} is immediate from \eqref{eq:j} and \eqref{eq:GH}. The function $G_\varphi$ is convex since $\delta_\varphi\ge0$ by Lemma~\ref{lem:varphi_prime_nneg_ninc}. For $H_\varphi$, it suffices that each block $x\mapsto j_\varphi(\|x\|)$ is convex on $\mathbb{R}^2$, which holds when $j_\varphi$ is convex and nondecreasing on $[0,\infty)$, as the composition of a convex nondecreasing function with the norm. Convexity of $j_\varphi$ follows from \eqref{eq:j}, since $\varphi$ is concave. For $0<s<t$, the secant slopes of the concave function $\varphi$ are nonincreasing, and together with $\varphi(0)=0$ and $\varphi(s)\le\delta_\varphi s$ (Lemma~\ref{lem:varphi_prime_nneg_ninc}) this gives
\[
  \frac{\varphi(t)-\varphi(s)}{t-s}\;\le\;\frac{\varphi(s)-\varphi(0)}{s-0}\;\le\;\delta_\varphi ;
\]
for $s=0$, the same bound is $\varphi(t)\le\delta_\varphi t$. Therefore,
$j_\varphi(t)-j_\varphi(s)=\delta_\varphi(t-s)-\bigl(\varphi(t)-\varphi(s)\bigr)\ge0$ for
$0\le s\le t$, i.e., $j_\varphi$ is nondecreasing.
\end{proof}

Without smoothing, the construction of Proposition~\ref{prop:sufficiency} fails for every weight of the total variation term.

\begin{cor}
\label{cor:no-splitting}
Let $q\in(0,1)$. For every $\lambda\ge0$, the function
\[
  H(u):=\lambda\sum_{i=1}^m\phi_i(u)-\beta\sum_{i=1}^m\phi_i(u)^q
\]
is not convex on $\mathbb{R}^n$. Consequently, the objective of \eqref{eq:tvq_problem} admits no decomposition $J=G-H$ with $G$, $H$ convex and $G=\tfrac12\|A\cdot-f\|^2+\lambda\sum_{i=1}^m\phi_i$.
\end{cor}

\begin{proof}
  We first fix $e \in \mathbb{R}^n$ with $\mathbb{K}e \neq 0$ and set $a_i := \phi_i(e)$, so that $a_i\geq 0$ and $\sum_{i=1}^{m}a_i^q>0$. Next, we evaluate the function at the symmetric points $se$ and $-se$ for a scalar $s > 0$. Then $\phi_i(\pm se) = s a_i$ and substituting this into the function we obtain
  \begin{equation*}
    H(se) = H(-se) = s\lambda \sum_{i=1}^m a_i - \beta s^q \sum_{i=1}^m a_i^q.
  \end{equation*}

  Given that the exponent $q \in (0,1)$, for a sufficiently small value of $s$, the fractional power $s^q$ is significantly larger than $s$. Therefore 
  $$
  H(se)=H(-se)<0=H(0)
  $$
  If we were to assume that $H$ is convex, by definition its value at the midpoint should satisfy the inequality 
  $$H(0) \le \frac{1}{2}H(se) + \frac{1}{2}H(-se)<0=H(0),$$
  which is a contradiction.

  Consequently, if $J = G - H$ with $G$ as stated, then algebraically $H = G - J$, which is exactly the function $H$ analyzed above. Since it has already been shown that $H$ is not convex, we conclude that such a DC decomposition is impossible.
\end{proof}
Finally, the component $H_\varphi$ of \eqref{eq:GH} is continuously differentiable, with a gradient bounded blockwise by the weight $\lambda=\beta\delta_\varphi$.

\begin{lem}
\label{lem:H_dif}
Let $\varphi$ be admissible. Then $j_\varphi\in C^1([0,\infty))$, with $j_\varphi'=\delta_\varphi-\varphi'$ on $(0,\infty)$, $j_\varphi'(0)=0$ and $0\le j_\varphi'\le\delta_\varphi$ on $[0,\infty)$. Consequently, $H_\varphi$ is continuously differentiable on $\mathbb{R}^n$, with
\begin{equation}
  \nabla H_\varphi(u)=\mathbb{K}^\top w_\varphi(u),
  \qquad
  w_\varphi(u)_i:=\beta\,j_\varphi'\bigl(\phi_i(u)\bigr)\,
  \frac{(\mathbb{K}u)_i}{\phi_i(u)}\in\mathbb{R}^2,
  \quad i=1,\dots,m,
  \label{w}
\end{equation}
where $w_\varphi(u)_i:=0$ when $\phi_i(u)=0$, and
\begin{equation}
  \|w_\varphi(u)_i\|=\beta\,j_\varphi'\bigl(\phi_i(u)\bigr)\;\le\;\beta\delta_\varphi=\lambda .
  \label{eq:tilt-bound}
\end{equation}
\end{lem}

\begin{proof}
The derivative of $j_\varphi$ at the origin follows from \ref{ass:A3} and Lemma~\ref{lem:varphi_prime_nneg_ninc}. The derivative of each block then follows from the chain rule, and at the degenerate blocks with $\phi_i(u)=0$ the estimate $j_\varphi(s)-j_\varphi(0)=o(s)$ gives the derivative $0$ directly; see Appendix~\ref{app:proofs-splitting}.
\end{proof}

The block $w_\varphi(u)_i$ has the direction of the gradient at pixel $i$, and its magnitude never exceeds the weight of the total variation term in $G_\varphi$. It is the vector that the algorithm of Section~\ref{sec:pdhg} subtracts from each subproblem.

We can now state the optimality system of \eqref{eq:Pphi} for every admissible penalty.

\begin{teo}\label{teo:opt_cond}
Let $\varphi$ be admissible and $\lambda=\beta\delta_\varphi$. If $u^*$ is a local minimizer of \eqref{eq:Pphi}, then there exist $z\in\mathbb{R}^{2m}$ such that $u^*$ satisfies the following optimality system:
\begin{subequations}\label{eq:sis_opt_low}
\begin{align}
    \mathbb{K}^\top w_\varphi(u^*)-A^\top(Au^*-f)-\lambda\mathbb{K}^\top z=0, \label{sist1}\\
    z_i\in \partial\|\cdot\|\bigl((\mathbb{K}u^*)_i\bigr),\quad i=1,\ldots,m,\label{sist2}
\end{align}
\end{subequations}
where $z_i\in\mathbb{R}^2$ denotes the $i$-th block of $z$, ordered as the blocks of $w_\varphi$ in \eqref{w}.
\end{teo}
\begin{proof}
See Appendix~\ref{app:proofs-splitting}.
\end{proof}

Since $H_\varphi$ is differentiable, \eqref{eq:sis_opt_low} is equivalent to $\partial G_\varphi(u^*)\cap\partial H_\varphi(u^*)\neq\emptyset$, i.e., to $u^*$ being a critical point of $G_\varphi-H_\varphi$ in the sense of difference-of-convex programming~\cite{hor}.

For the Huber-type smoothing, the system \eqref{eq:sis_opt_low} is written out block by block in Corollary~\ref{cor:opt_huber} below.

\subsection{DC representation of the Huber-type smoothing}
\label{sec:huber_smoothing}
The previous results hold for any admissible penalty. In this subsection, we apply them to the Huber-type smoothing $\varphi_\gamma^q$ of \eqref{eq:reg}, which is admissible by Lemma~\ref{lem:huber-admissible}. Consequently, all results of this section apply to it. To simplify the notation when using $\varphi_\gamma^q$, we use the abbreviations
\begin{equation*}
  J_\gamma:=J_{\varphi_\gamma^q},\quad
  G_\gamma:=G_{\varphi_\gamma^q},\quad
  H_\gamma:=H_{\varphi_\gamma^q},\quad
  j_\gamma:=j_{\varphi_\gamma^q},\quad
  w:=w_{\varphi_\gamma^q},\quad
  \lambda=\beta\delta_\varphi ,
\end{equation*}
so that $J_\gamma=G_\gamma-H_\gamma$ is the splitting \eqref{eq:dc} for this penalty. The function $j_\gamma$ of \eqref{eq:j}, from which $H_{\gamma}$ is built and which the algorithm linearizes, is explicit:

\begin{equation}
  j_\gamma(t):=\delta_\varphi t-\varphi_\gamma^q(t)
  =\left\{\begin{array}{ll}
     \delta_{\varphi}t-\left(t+\frac{q-1}{\gamma}\right)^{q}, & t>1\slash\gamma,\\[6pt]
     0, & t\leq 1\slash\gamma.
   \end{array}\right.
  \label{j}
\end{equation}
Being the gap between the concave $\varphi_\gamma^q$ and its tangent at the origin, $j_\gamma$ is nonnegative and convex, and vanishes up to the knot; the concavity of the penalty sits in $-j_\gamma$. By Lemma~\ref{lem:H_dif} it is $C^1$ and nondecreasing, with
\begin{equation}
  j_\gamma'(t)=\left\{\begin{array}{ll}
     \delta_{\varphi}-q\left(t+\frac{q-1}{\gamma}\right)^{q-1}, & t>1\slash\gamma,\\[6pt]
     0, & t\leq 1\slash\gamma,
   \end{array}\right.
  \label{jprima}
\end{equation}
and $0\le j_\gamma'\le\delta_\varphi$ on $[0,\infty)$.

By Lemma~\ref{lem:H_dif}, $H_\gamma$ is continuously differentiable on $\mathbb{R}^n$ with $\nabla H_\gamma(u)=\mathbb{K}^\top w(u)$, where $w(u)$ is \eqref{w} with $j_\gamma'$ as in \eqref{jprima}. Since $j_\gamma'\equiv0$ on $[0,1/\gamma]$, the block $w(u)_i$ vanishes whenever $\phi_i(u)\le1/\gamma$, and $\|w(u)_i\|\le\lambda$ otherwise by \eqref{eq:tilt-bound}. Theorem~\ref{teo:opt_cond} therefore applies to \eqref{eq:Pgamma} with $w_\varphi=w$ and $\lambda=\beta q^q\gamma^{1-q}$, and substituting \eqref{jprima} into \eqref{w} gives its optimality system in explicit form.

\begin{cor}
\label{cor:opt_huber}
Let $q\in(0,1)$, $\gamma\gg1$ and $\lambda=\beta q^q\gamma^{1-q}$. If $u^*$ is a local
minimizer of \eqref{eq:Pgamma}, then there exists $z\in\mathbb{R}^{2m}$ such that
\begin{subequations}\label{eq:sis_opt_huber}
\begin{align}
  A^\top(Au^*-f)+\lambda\,\mathbb{K}^\top z-\mathbb{K}^\top w(u^*)&=0,
  \label{eq:huber_sist1}\\
  z_i\in\partial\|\cdot\|\bigl((\mathbb{K}u^*)_i\bigr),\quad i=1,\dots,m,&
  \label{eq:huber_sist2}
\end{align}
\end{subequations}
where the blocks of $w(u^*)$ are
\begin{equation}
  w(u^*)_i=\left\{\begin{array}{ll}
    \beta\Bigl(q^q\gamma^{1-q}-q\Bigl(\phi_i(u^*)+\dfrac{q-1}{\gamma}\Bigr)^{q-1}\Bigr)
    \dfrac{(\mathbb{K}u^*)_i}{\phi_i(u^*)}, & \phi_i(u^*)>1/\gamma,\\[12pt]
    0, & \phi_i(u^*)\le1/\gamma.
  \end{array}\right.
  \label{eq:w_huber}
\end{equation}
Equivalently, $A^\top(Au^*-f)+\mathbb{K}^\top\zeta=0$ for some $\zeta\in\mathbb{R}^{2m}$
with blocks
\begin{equation}
  \zeta_i=\left\{\begin{array}{ll}
    \beta\,q\Bigl(\phi_i(u^*)+\dfrac{q-1}{\gamma}\Bigr)^{q-1}
    \dfrac{(\mathbb{K}u^*)_i}{\phi_i(u^*)}, & \phi_i(u^*)>1/\gamma,\\[12pt]
    \lambda\,\dfrac{(\mathbb{K}u^*)_i}{\phi_i(u^*)}, & 0<\phi_i(u^*)\le1/\gamma,\\[12pt]
    \lambda\,z_i,\quad \|z_i\|\le1, & \phi_i(u^*)=0.
  \end{array}\right.
  \label{eq:zeta_huber}
\end{equation}
\end{cor}

\begin{proof}
The penalty $\varphi_\gamma^q$ is admissible with $\delta_\varphi=q^q\gamma^{1-q}$ (Lemma~\ref{lem:huber-admissible}), so Theorem~\ref{teo:opt_cond} gives \eqref{eq:sis_opt_huber}, and \eqref{eq:w_huber} is \eqref{w} with $j_\gamma'$ as in \eqref{jprima}. For the second form, set $\zeta_i:=\lambda z_i-w(u^*)_i$. If $\phi_i(u^*)>0$, then $\partial\|\cdot\|\bigl((\mathbb{K}u^*)_i\bigr)$ is the singleton $\{(\mathbb{K}u^*)_i/\phi_i(u^*)\}$, and $\lambda-\beta j_\gamma'(\phi_i(u^*))=\beta(\varphi_\gamma^q)'(\phi_i(u^*))$, which equals $\lambda$ below the knot and $\beta q\bigl(\phi_i(u^*)+\frac{q-1}{\gamma}\bigr)^{q-1}$ above it by \eqref{eq:reg-deriv}. If $\phi_i(u^*)=0$, then $w(u^*)_i=0$ and $\partial\|\cdot\|(0)$ is the closed unit ball of $\mathbb{R}^2$.
\end{proof}

On the blocks with $\phi_i(u^*)>0$, $\mathbb{K}_i^\top\zeta_i$ is the gradient of the continuously differentiable map $u\mapsto\beta\varphi_\gamma^q(\phi_i(u))$ at $u^*$, so \eqref{eq:sis_opt_huber} is the classical first-order condition there; the first branch of \eqref{eq:zeta_huber} is the gradient of $\beta\phi_i(\cdot)^q$ with the argument shifted by $(q-1)/\gamma$. The nonsmoothness of \eqref{eq:Pgamma} is confined to the blocks with $\phi_i(u^*)=0$, where the penalty coincides with $\lambda\|\cdot\|$ near $u^*$ and the multiplier $\zeta_i$ ranges over the ball of radius $\lambda$, exactly as in the total variation problem with weight $\lambda$. The system \eqref{eq:sis_opt_huber} states that $u^*$ is a critical point of $G_\gamma-H_\gamma$. This is the stationarity notion of the algorithm of Section~\ref{sec:pdhg}, up to the multiplier of the box that the algorithm adds to its subproblems (Corollary~\ref{cor:opt_box} and Theorem~\ref{teo:convergence}).
 
\section{Smoothing error analysis}
\label{sec:smoothing_error}
The algorithm of Section~\ref{sec:pdhg} minimizes $J_\gamma$ instead of the objective of \eqref{eq:tvq_problem}. The error of this replacement decreases like $\gamma^{-q}$, so a minimizer of $J_\gamma$ is an approximate minimizer of the original problem (Lemma~\ref{lem:smoothing-error}). The subproblem weight $\lambda=\beta\delta_\varphi$, however, grows like $\gamma^{1-q}$. The following lemma quantifies the error of the smoothing $\varphi_\gamma^q$ of \eqref{eq:reg} and its effect on the objective.

\begin{lem}
\label{lem:smoothing-error}
For every $t\ge0$, $q\in(0,1)$ and $\gamma>0$,
\begin{equation}
  0\;\le\;t^q-\varphi_\gamma^q(t)\;\le\;q^q(1-q)\,\gamma^{-q},
  \label{eq:smoothing-error}
\end{equation}
and the upper bound is attained, at $t=q/\gamma$. Consequently, writing
$J(u)=\tfrac12\|Au-f\|^2+\beta\sum_{i=1}^m\|(\mathbb{K}u)_i\|^q$ for the unsmoothed
objective,
\begin{equation}
  0\;\le\;J(u)-J_\gamma(u)\;\le\;E(\gamma):=\beta\,m\,q^q(1-q)\,\gamma^{-q}
  \qquad\text{for all }u\in \R^n,
  \label{eq:E-gamma}
\end{equation}
and any minimizer $u_\gamma$ of $J_\gamma$ satisfies $J(u_\gamma)\le J(u)+E(\gamma)$ for every $u\in \R^n$, i.e., $J(u_\gamma)\le\inf J+E(\gamma)$.
\end{lem}

\begin{proof}
See Appendix~\ref{app:proofs-splitting}.
\end{proof}

\begin{rem}
No admissible penalty with the same slope is more accurate than $\varphi_\gamma^q$. Since $\varphi\le\delta_\varphi t$ for every admissible $\varphi$ (Lemma~\ref{lem:varphi_prime_nneg_ninc}), its error $t^q-\varphi(t)$ is at least $t^q-\delta_\varphi t$, whose maximum over $t\ge0$ is $(1-q)q^{q/(1-q)}\delta_\varphi^{-q/(1-q)}$. For $\delta_\varphi=q^q\gamma^{1-q}$ this is the bound $q^q(1-q)\gamma^{-q}$ of Lemma~\ref{lem:smoothing-error}, which $\varphi_\gamma^q$ attains at $t=q/\gamma$.
\end{rem}

\section{The Primal-Dual DC  method}
\label{sec:pdhg}
In this section, we propose an algorithm for the numerical solution of \eqref{eq:Pgamma} based on the DC algorithm (DCA). We restrict the problem to the dynamic range of the image,
\begin{equation}
  B:=\{u\in\R^n:\,0\le u_i\le 1,\,i=1,\dots,n\},
  \label{eq:box}
\end{equation}
whose indicator function we denote by $\iota_B$. The box is not part of~\eqref{eq:Pgamma}, but it contains every image the model is meant to recover. We add it so that the duality certificate of Proposition~\ref{prop:gap} is finite for every operator $A$; see Remark~\ref{rem:no-certificate}.

The DCA is based on the DC criticality condition $\partial H_\gamma (u^*) \subset \partial (G_\gamma + \iota_B)(u^*) $ which induces the following primal--dual scheme \cite{taole}:
  
\begin{align*}
	& \text{first choose: } v_k \in  \partial H_\gamma (u_k), \\
	&  \text{then choose: } u_{k+1} \in  \partial (G_\gamma + \iota_B)^*(v_k),
\end{align*}
where $\partial H_\gamma$ denotes the subdifferential of $H_\gamma$ and $ \partial (G_\gamma + \iota_B)^*$ corresponds to the subdifferential of the Fenchel conjugate of $G_\gamma + \iota_B$. The DCA therefore replaces, at each iteration, the concave part $-H_\gamma$ by its linearization at the current iterate \cite{hor}. Since $H_\gamma$ is differentiable by Lemma~\ref{lem:H_dif}, this linearization only requires the evaluation of $v_k:=\mathbb{K}^\top w(u_k)$, and the outer step reads
\begin{equation}
  u_{k+1}\in\argmin_{u\in B}\;
  \frac12\|Au-f\|^2+\lambda\sum_{i=1}^m \phi_i(u)-\scalar{v_k}{u},
  \qquad \lambda=\beta\delta_\varphi.
  \label{eq:subproblem}
\end{equation}
Since $B$ is compact, the convex problem~\eqref{eq:subproblem} has a minimizer.

Algorithm~\ref{alg:dca} therefore approximates critical points of $J_\gamma$ over $B$ rather than of \eqref{eq:Pgamma}. Their optimality system differs from \eqref{eq:sis_opt_huber} only by the multiplier of the box. The following corollary of Theorem~\ref{teo:opt_cond} states the optimality system of a local minimizer of $J_\gamma$ over $B$.

\begin{cor}
\label{cor:opt_box}
If $u^*$ is a local minimizer of $J_\gamma$ over $B$, then there exist
$z\in\mathbb{R}^{2m}$ and $\nu\in N_B(u^*)$ such that
\begin{equation}
  \mathbb{K}^\top w(u^*)-A^\top(Au^*-f)-\lambda\mathbb{K}^\top z-\nu=0,\qquad
  z_i\in\partial\|\cdot\|\bigl((\mathbb{K}u^*)_i\bigr),\quad i=1,\dots,m,
  \label{eq:sis_opt_box}
\end{equation}
where $N_B(u):=\{\nu\in\mathbb{R}^n:\scalar{\nu}{y-u}\le0\ \text{for all }y\in B\}$ is the normal cone to $B$ at $u$, i.e., $\nu_i\le0$ where $u_i=0$, $\nu_i\ge0$ where $u_i=1$, and $\nu_i=0$ otherwise. If $\nu=0$, in particular if $u^*$ lies in the interior of $B$, \eqref{eq:sis_opt_box} is the optimality system \eqref{eq:sis_opt_huber} of \eqref{eq:Pgamma}.
\end{cor}

\begin{proof}
The proof of Theorem~\ref{teo:opt_cond} applies to $G_\gamma+\iota_B-H_\gamma$. Since $G_\gamma$ is finite and continuous, the sum rule adds $\partial\iota_B(u^*)=N_B(u^*)$ to $\partial G_\gamma(u^*)$ \cite{urru}.
\end{proof}

The general algorithm is presented below:
\begin{algorithm}[h]
\caption{Primal-Dual DC algorithm for \eqref{eq:Pgamma} over the box $B$}
\label{alg:dca}
\begin{algorithmic}[1]
\Require $u_0\in B$, $\lambda=\beta\delta_\varphi$, tolerance $\varepsilon\ge0$, ratio
         $\rho\in(0,1)$, balance $c_0>0$
\For{$k=0,1,2,\ldots$}
  \State $v_k \gets \mathbb{K}^\top w(u_k)$
  \State $(\hat u_{k+1},r_{k+1},s_{k+1}) \gets$ Algorithm~\ref{alg:pdhgm}$(v_k,u_k,\rho,\varepsilon)$ with $\tau=c_k\bar\tau$, $\sigma=\bar\tau/c_k$
  \State $c_{k+1} \gets$ {\scriptsize the balance \eqref{eq:balance}}
  \If{$\eta_k:=\mathcal{G}_k(u_k,r_{k+1},s_{k+1})\le\varepsilon$}
    \State \Return $u_k$
  \EndIf
  \State $u_{k+1} \gets \hat{u}_{k+1}$
\EndFor
\end{algorithmic}
\end{algorithm}

All subproblems in Algorithm \ref{alg:dca} share the total variation term with the scalar weight $\lambda$, which is fixed before the run, and differ only in the linear term $-\scalar{v_k}{u}$. One bound on the step sizes therefore serves all of them, and a balance between the primal and dual steps that suits one subproblem also suits the next. By \eqref{eq:tilt-bound}, the blocks of $w(u_k)$ are bounded by $\lambda$, i.e., $w(u_k)$ lies in the dual ball of the total variation term. Algorithm~\ref{alg:dca} implements the outer loop, and Algorithm~\ref{alg:pdhgm} in Section~\ref{sec:saddle} solves its subproblems.

\subsection{Saddle-point formulation of the DC subproblem}
\label{sec:saddle}
If the fidelity term is kept with the primal variable, every inner iteration requires $(I+\tau A^\top A)^{-1}$, whose computation is expensive or unavailable when $A$ is a convolution, a coordinate projection or an X-ray transform. We avoid this inversion by dualizing $A$ and $\mathbb{K}$ together. To this end, we define the operator $L$ and the functions $\tilde F$ and $\tilde G$ by:
\begin{equation*}
  L=\begin{bmatrix} A\\ \mathbb{K}\end{bmatrix},\qquad
  \tilde F(r,s)=\underbrace{\tfrac12\|r-f\|^2}_{\tilde F_1(r)}
               +\underbrace{\lambda\sum_{i=1}^m\|s_i\|}_{\tilde F_2(s)},\qquad
  \tilde G(u)=-\scalar{v_k}{u}+\iota_B(u),
\end{equation*}
where $r\in\mathbb{R}^M$ and $s\in\mathbb{R}^{2m}$. Therefore, problem \eqref{eq:subproblem} is expressed by
\begin{equation*}
  \min_u \tilde F(Lu)+\tilde G(u).
\end{equation*}
Since $\tilde F$ is proper, convex and lower semicontinuous, it coincides with its biconjugate, $\tilde F(Lu)=\max_{r,s}\,\scalar{Lu}{(r,s)}-\tilde F^*(r,s)$, where the conjugate can be expressed as $\tilde F^*(r,s)=\tfrac12\|r\|^2+\scalar{f}{r}+\sum_{i=1}^m\iota_{\{\|\cdot\|\le\lambda\}}(s_i)$. Substituting this representation of $\tilde F(Lu)$ yields the saddle-point problem:
\begin{equation}
  \min_{u\in B}\ \max_{r\in\mathbb{R}^M,\,s\in\mathbb{R}^{2m}}\;
  \scalar{Au}{r}+\scalar{\mathbb{K}u}{s}-\scalar{v_k}{u}
  -\tfrac12\|r\|^2-\scalar{f}{r}
  -\sum_{i=1}^m\iota_{\{\|\cdot\|\le\lambda\}}(s_i).
  \label{eq:saddle_point}
\end{equation}
This problem has the structure considered by Chambolle and Pock \cite{chambolle2016}, whose primal--dual method only requires the proximal maps of $\tilde F^*$ and $\tilde G$. Throughout, we write
\begin{equation}
  \mathcal{S}:=\bigl\{s\in\mathbb{R}^{2m}\;:\;\|s_i\|\le\lambda,\ i=1,\dots,m\bigr\}
  \label{eq:dualball}
\end{equation}
for the dual ball of the total variation term, so that $\tilde F_2^*=\iota_{\mathcal{S}}$ and the corresponding dual step is the projection $\proj_{\mathcal{S}}$.

A direct advantage in this formulation is that all proximal maps are available in closed form. Indeed, we have
\begin{equation*}
  \prox_{\sigma\tilde F_1^*}(r)=\frac{r-\sigma f}{1+\sigma},
  \qquad
  \bigl(\prox_{\sigma\tilde F_2^*}(s)\bigr)_i=\frac{\lambda\,s_i}{\max(\lambda,\|s_i\|)},
  \qquad
  \prox_{\tau\tilde G}(u)=\proj_B(u+\tau v_k) .
\end{equation*}
The primal proximal step reduces to a translation followed by a projection onto $B$, so $A$ and $A^\top$ are only applied.

The step sizes satisfy the usual condition $\tau\sigma\|L\|^2<1$, where $\|L\|^2\le\|A\|^2+\|\mathbb{K}\|^2$ and $\|\mathbb{K}\|^2\le8$ for forward differences. We take $\tau=c\bar\tau$ and $\sigma=\bar\tau/c$ with $\bar\tau=0.99/\sqrt{\|A\|^2+8}$, which satisfy the condition for every balance $c>0$. Since $L$ does not change between outer iterations, $\bar\tau$ is computed once and used for all subproblems.

The balance $c$ does not affect the convergence of the method, but its speed depends strongly on it \cite[p.~54]{chambolle2016}. We estimate a suitable balance by the displacement over the previous subproblem and set
\begin{equation}
  c_{k+1}:=\Bigl(c_k\,\frac{\|\hat u_{k+1}-u_k\|}{\|(r_{k+1},s_{k+1})-(r_k,s_k)\|}\Bigr)^{1/2},
  \label{eq:balance}
\end{equation}
i.e., the geometric mean of the old balance and the estimate. We keep $c_{k+1}$ in $[10^{-2},10^{2}]$ and set $c_{k+1}=c_k$ if either displacement vanishes. The first balance $c_0$ is an input of Algorithm~\ref{alg:dca}; its choice in the experiments is described in Section~\ref{sec:setup}. Since the balance is fixed within each subproblem, none of the results below depends on $c_k$.

Using the proximal operators above, the subproblem \eqref{eq:subproblem} is solved by the following algorithm. 

\begin{algorithm}[h]
\caption{Algorithm for the subproblem \eqref{eq:subproblem}}
\label{alg:pdhgm}
\begin{algorithmic}[1]
\Require tilt $v_k$, current iterate $u_k\in B$, ratio $\rho\in(0,1)$, tolerance
         $\varepsilon\ge0$, steps $\tau,\sigma>0$ with $\tau\sigma\|L\|^2<1$,
         $(u^0,r^0,s^0)$, $\bar u^0=u^0$
\State $p_k \gets \Phi_k(u_k)$
\For{$l=0,1,2,\ldots$}
  \State $r^{l+1} \gets \bigl(r^l+\sigma(A\bar u^l-f)\bigr)/(1+\sigma)$
  \State $s^{l+1} \gets \proj_{\mathcal{S}}\bigl(s^l+\sigma\mathbb{K}\bar u^l\bigr)$
  \State $u^{l+1} \gets \proj_B\bigl(u^l-\tau(A^\top r^{l+1}+\mathbb{K}^\top s^{l+1}-v_k)\bigr)$
  \State $d^{l+1} \gets D_k(r^{l+1},s^{l+1})$
  \If{$p_k-\Phi_k(u^{l+1})\ge(1-\rho)\,(p_k-d^{l+1})$ \textbf{or}
      $p_k-d^{l+1}\le\varepsilon$}
    \State \Return $(u^{l+1},r^{l+1},s^{l+1})$
  \EndIf
  \State $\bar u^{l+1} \gets 2u^{l+1}-u^l$
\EndFor
\end{algorithmic}
\end{algorithm}
Algorithm~\ref{alg:pdhgm} is started from $(u^0,r^0,s^0)=(u_k,r_k,s_k)$, i.e., from the current iterate and the dual variables returned by the previous subproblem, with $(r_0,s_0)=(0,0)$.

Next, we provide a stopping criterion of Algorithm~\ref{alg:pdhgm}. The current iterate $u_k$ is a critical point exactly when the convex model of $J_\gamma$ built at $u_k$ cannot be decreased any further (Corollary~\ref{cor:opt_box}), so the quantity of interest is the subproblem error at the
\emph{current} iterate,
\begin{equation}
  \theta_k:=\Phi_k(u_k)-\min\Phi_k\;\ge\;0 .
  \label{eq:theta}
\end{equation}

Since $\Phi_k=G_\gamma+\iota_B-\scalar{v_k}{\cdot}$ is convex, from \eqref{eq:theta} we have that  $\theta_k=0$ if and only if $u_k$ minimizes $\Phi_k$, i.e., $v_k\in\partial(G_\gamma+\iota_B)(u_k)$.

Moreover, $G_\gamma$ is finite and continuous, the sum rule gives $\partial(G_\gamma+\iota_B)(u_k)=\partial G_\gamma(u_k)+N_B(u_k)$ \cite{urru}. With the expression of $\partial G_\gamma(u_k)$ given in the proof of Theorem~\ref{teo:opt_cond} and $v_k=\nabla H_\gamma(u_k)=\mathbb{K}^\top w(u_k)$, this inclusion is \eqref{eq:sis_opt_box} at $u_k$. Thus $\theta_k$ measures stationarity. However, this error cannot be computed, since $\min\Phi_k$ is unknown. The next proposition bounds $\theta_k$ by a computable certificate, which vanishes at every saddle point of~\eqref{eq:saddle_point}.

\begin{prop}
\label{prop:gap}
Let $\Phi_k$ denote the objective of \eqref{eq:subproblem},
\begin{equation}
  \Phi_k(u):=\tfrac12\|Au-f\|^2+\lambda\sum_{i=1}^m \phi_i(u)-\scalar{v_k}{u}+\iota_B(u),
  \qquad v_k=\mathbb{K}^\top w(u_k),
  \label{eq:Phik}
\end{equation}
and for $r\in\mathbb{R}^M$ and $s\in\mathcal{S}$ define the dual value
\begin{equation}
  D_k(r,s):=-\tfrac12\|r\|^2-\scalar{f}{r}-\sigma_B(\xi),
  \qquad \xi:=v_k-A^\top r-\mathbb{K}^\top s ,
  \label{eq:dualvalue}
\end{equation}
where $\sigma_B(\xi):=\sup_{u\in B}\scalar{\xi}{u}=\sum_{i=1}^n\max(\xi_i,0)$ is the support function of the box~\eqref{eq:box}. Then $D_k(r,s)\le\min\Phi_k$, and for every $u\in B$,
\begin{equation}
  0\;\le\;\Phi_k(u)-\min\Phi_k\;\le\;\mathcal{G}_k(u,r,s):=\Phi_k(u)-D_k(r,s).
  \label{eq:gap}
\end{equation}
Moreover, $\mathcal{G}_k$ vanishes at every saddle point of~\eqref{eq:saddle_point}.
\end{prop}

\begin{proof}
Let $\mathcal{L}(u,r,s)=\scalar{Au}{r}+\scalar{\mathbb{K}u}{s}-\scalar{v_k}{u}
  -\tfrac12\|r\|^2-\scalar{f}{r}
  -\sum_{i=1}^m\iota_{\{\|\cdot\|\le\lambda\}}(s_i)$ denote the function in \eqref{eq:saddle_point}, for $u\in B$, $r\in\mathbb{R}^M$ and $s\in\mathcal{S}$. Both claims follow from two partial optimizations of $\mathcal{L}$.

First, fix $u\in B$. The terms in $r$ form the concave quadratic $\scalar{Au-f}{r}-\tfrac12\|r\|^2$, which is maximized at $r=Au-f$ with value $\tfrac12\|Au-f\|^2$. By the Cauchy--Schwarz inequality, $\scalar{(\mathbb{K}u)_i}{s_i}\le\lambda\|(\mathbb{K}u)_i\|=\lambda\phi_i(u)$ whenever $\|s_i\|\le\lambda$, with equality for $s_i=\lambda(\mathbb{K}u)_i/\|(\mathbb{K}u)_i\|$ (any $s_i$ with $\|s_i\|\le\lambda$ if $(\mathbb{K}u)_i=0$). Hence
\begin{equation*}
  \max_{r\in\mathbb{R}^M,\,s\in\mathcal{S}}\mathcal{L}(u,r,s)=\Phi_k(u).
\end{equation*}
On the other hand, fix $r\in\mathbb{R}^M$ and $s\in\mathcal{S}$. Collecting the terms in $u$ gives $\mathcal{L}(u,r,s)=-\scalar{\xi}{u}-\tfrac12\|r\|^2-\scalar{f}{r}$ with $\xi$ given by \eqref{eq:dualvalue}, and since the minimum of $-\scalar{\xi}{u}$ over $B$ is $-\sigma_B(\xi)$,
\begin{equation*}
  \min_{u\in B}\mathcal{L}(u,r,s)=D_k(r,s).
\end{equation*}
Together, these identities give $\Phi_k(u)\ge\mathcal{L}(u,r,s)\ge D_k(r,s)$ for every $u\in B$. Taking the minimum over $u$ yields $D_k(r,s)\le\min\Phi_k$, and subtracting this inequality from $\Phi_k(u)$ yields \eqref{eq:gap}.

Let $(u^*,r^*,s^*)$ be a saddle point, i.e., $(r^*,s^*)$ maximizes $\mathcal{L}(u^*,\cdot,\cdot)$ and $u^*$ minimizes $\mathcal{L}(\cdot,r^*,s^*)$ over $B$. By this property and the two identities,
\begin{equation*}
  \Phi_k(u^*)=\max_{r\in\mathbb{R}^M,\,s\in\mathcal{S}}\mathcal{L}(u^*,r,s)=\mathcal{L}(u^*,r^*,s^*)=\min_{u\in B}\mathcal{L}(u,r^*,s^*)=D_k(r^*,s^*).
\end{equation*}
Therefore, $\mathcal{G}_k(u^*,r^*,s^*)=0$. A saddle point exists by the proof of Lemma~\ref{lem:inner-stop}.
\end{proof}

Every dual pair $(r,s)$ therefore gives a lower bound on $\min\Phi_k$, and $\mathcal{G}_k$ bounds how far $\Phi_k(u)$ is from the optimal value. The bound is finite for every $(r,s)$, and its evaluation costs one application of $A^\top$ and one of $\mathbb{K}^\top$ plus $O(n)$ operations. All iterates of Algorithm~\ref{alg:pdhgm} are admissible dual pairs, since the projection places $s$ in $\mathcal{S}$ and $r$ is unconstrained.

\begin{rem}
\label{rem:no-certificate}
Notice that the certificate is finite because $B$ is bounded. Without the box, the minimum of $-\scalar{\xi}{u}$ in the proof of Proposition~\ref{prop:gap} is $-\infty$ unless $\xi=0$, i.e., unless $(r,s)$ lies in the affine set $\{A^\top r+\mathbb{K}^\top s=v_k\}$, which the iterates of Algorithm~\ref{alg:pdhgm} reach only in the limit. The dual value is then $-\infty$ at every iterate, and \eqref{eq:gap} gives no information. Correcting $r$ to restore feasibility requires $A$ to be injective, which excludes a coordinate projection and an X-ray transform. Correcting $s$ requires a Poisson solve and a rescaling into $\mathcal{S}$ that weakens the bound when $\beta$ is small. For a bounded $B$, $\sigma_B$ is finite for every $A$, at the cost of one clip per inner iteration.
\end{rem}

 Given $\rho\in(0,1)$ and the tolerance $\varepsilon\ge0$ of Algorithm~\ref{alg:dca}, Algorithm~\ref{alg:pdhgm} stops at the first triple $(u^{l+1},r^{l+1},s^{l+1})$ that satisfies
\begin{equation}
  \Phi_k(u_k)-\Phi_k(u^{l+1})\;\ge\;(1-\rho)\,\mathcal{G}_k(u_k,r^{l+1},s^{l+1})
  \quad\text{or}\quad
  \mathcal{G}_k(u_k,r^{l+1},s^{l+1})\;\le\;\varepsilon .
  \label{eq:relstop}
\end{equation}
The certificate $\mathcal{G}_k(u_k,r^{l+1},s^{l+1})$ bounds the decrease that the subproblem can still achieve from $u_k$, and the first test accepts $u^{l+1}$ once it attains a fraction $1-\rho$ of this bound. Since the bound shrinks as $u_k$ approaches stationarity, the test tightens along the iteration and no sequence of tolerances has to be chosen in advance. The second test stops the inner loop once the certificate shows that $\theta_k\le\varepsilon$. As the certificates at $u_k$ and $u^{l+1}$ share the dual value, the first test can also be written as
\begin{equation*}
  \mathcal{G}_k(u^{l+1},r^{l+1},s^{l+1})\;\le\;\rho\,\mathcal{G}_k(u_k,r^{l+1},s^{l+1}),
\end{equation*}
which is the multiplicative inexactness criterion of \cite[eq.~(3)]{lee2019inexact} with $\min\Phi_k$ replaced by the dual value, as in \cite[eq.~(58)]{bonettini2016variable}. The next lemma shows that \eqref{eq:relstop} is reached after finitely many inner iterations, unless $u_k$ is already a critical point and $\varepsilon=0$.

\begin{lem}
\label{lem:inner-stop}
Let $u_k\in B$, $\rho\in(0,1)$ and $\varepsilon\ge0$, and assume that $\theta_k>0$ or $\varepsilon>0$. Then Algorithm~\ref{alg:pdhgm} returns after finitely many iterations a triple $(\hat u,r,s)\in B\times\mathbb{R}^M\times\mathcal{S}$ that satisfies \eqref{eq:relstop}. If it satisfies the first test in \eqref{eq:relstop}, then
\begin{equation}
  \Phi_k(u_k)-\Phi_k(\hat u)\;\ge\;(1-\rho)\,\mathcal{G}_k(u_k,r,s)\;\ge\;(1-\rho)\,\theta_k .
  \label{eq:model-decrease}
\end{equation}
\end{lem}

\begin{proof}
We first construct a saddle point of \eqref{eq:saddle_point}. Since $B$ is compact, \eqref{eq:subproblem} has a minimizer $u^*$. Therefore, applying the sum rule and since $G_\gamma$ is finite and continuous: $$v_k\in\partial(G_\gamma+\iota_B)(u^*) =\partial G_\gamma(u^*)+N_B(u^*).$$

By the expression of $\partial G_\gamma(u^*)$ in the proof of Theorem~\ref{teo:opt_cond}, $$v_k=A^\top(Au^*-f)+\lambda\mathbb{K}^\top z+\nu,$$ with $z_i\in\partial\|\cdot\|\bigl((\mathbb{K}u^*)_i\bigr)$ and $\nu\in N_B(u^*)$. Put $r^*=Au^*-f$ and $s^*=\lambda z\in\mathcal{S}$. For fixed $u^*$, the function in \eqref{eq:saddle_point} is maximized at $r^*$, where the gradient $Au^*-f-r$ of the concave quadratic vanishes, and at $s^*$, since $\scalar{(\mathbb{K}u^*)_i}{s_i}\le\lambda\|(\mathbb{K}u^*)_i\| =\scalar{(\mathbb{K}u^*)_i}{\lambda z_i}$ for $\|s_i\|\le\lambda$. For fixed $(r^*,s^*)$, it equals $-\scalar{\xi^*}{u}$ up to a constant, where $\xi^*=v_k-A^\top r^*-\mathbb{K}^\top s^*=\nu$ by \eqref{eq:dualvalue}, and it is minimized over $B$ at $u^*$ because $\nu\in N_B(u^*)$. Hence $(u^*,r^*,s^*)$ is a saddle point.

Since $\tau\sigma\|L\|^2<1$, the iterates of Algorithm~\ref{alg:pdhgm} converge to some saddle point, i.e., $B\times\mathbb{R}^M\times\mathcal{S}\ni(u^l,r^l,s^l)\to(\bar u,\bar r,\bar s)$ of \eqref{eq:saddle_point} as $l\to\infty$, see~\cite[Theorem 1]{cham}. 
Then, by Proposition~\ref{prop:gap}, $\mathcal{G}_k(\bar u,\bar r,\bar s)=0$, i.e., $D_k(\bar r,\bar s)=\Phi_k(\bar u)$, and since $D_k(\bar r,\bar s)\le\min\Phi_k\le\Phi_k(\bar u)$, both values equal $\min\Phi_k$. On the other hand, the function $D_k$ in \eqref{eq:dualvalue} is continuous, because $\sigma_B(\xi)=\sum_{i=1}^n\max(\xi_i,0)$ is, and $\mathcal{G}_k$ is continuous on $B\times\mathbb{R}^M\times\mathcal{S}$. Therefore, as $l\to\infty$,
\begin{equation*}
  \mathcal{G}_k(u^{l},r^{l},s^{l})\;\to\;0,
  \qquad
  \mathcal{G}_k(u_k,r^{l},s^{l})=\Phi_k(u_k)-D_k(r^{l},s^{l})\;\to\;\Phi_k(u_k)-\min\Phi_k=\theta_k .
\end{equation*}
Moreover, $\Phi_k(u^l)\to\Phi_k(\bar u)=\min\Phi_k$, since $\Phi_k$ is continuous on $B$.
If $\theta_k>0$, the left-hand side of the first inequality in \eqref{eq:relstop} tends to $\theta_k$ and its right-hand side to $(1-\rho)\,\theta_k<\theta_k$, so the first inequality is satisfied after finitely many iterations. If $\theta_k=0$ and $\varepsilon>0$, the left-hand side of the second inequality tends to $0<\varepsilon$, so the second test is passed after finitely many iterations.

Finally, the first test in \eqref{eq:relstop} is the first inequality in \eqref{eq:model-decrease}, and $\mathcal{G}_k(u_k,r,s)\ge\theta_k$ by Proposition~\ref{prop:gap}.
\end{proof}

The next proposition shows that the certificate $\eta_k$ of Algorithm~\ref{alg:dca} bounds the stationarity measure $\theta_k$ of \eqref{eq:theta}, and that $J_\gamma$ decreases by at least a fixed fraction of $\eta_k$ at every outer step. We refer to this property as \emph{certified descent}. It is the key ingredient for the convergence of Algorithm~\ref{alg:dca}.

\begin{prop}
\label{prop:descent}
Let $u_k\in B$, and let Algorithm~\ref{alg:pdhgm} return $(\hat u_{k+1},r_{k+1},s_{k+1})$ under the rule \eqref{eq:relstop}. Then $\theta_k\le\eta_k$. If, moreover, $\eta_k>\varepsilon$, i.e., if Algorithm~\ref{alg:dca} does not stop at iteration $k$, then $u_{k+1}=\hat u_{k+1}\in B$ and
\begin{align}
  \theta_k\;&\le\;\eta_k\;\le\;\frac{\theta_k}{1-\rho},
  \label{eq:eta-theta}\\
  J_\gamma(u_{k+1})\;&\le\;J_\gamma(u_k)-(1-\rho)\,\eta_k .
  \label{eq:descent}
\end{align}
\end{prop}

\begin{proof}
See Appendix~\ref{app:proofs-method}.
\end{proof}

With the descent of Proposition~\ref{prop:descent}, we can now prove the convergence of Algorithm~\ref{alg:dca}.

\begin{teo}
\label{teo:convergence}
Let $(u_k)$ be generated by Algorithm~\ref{alg:dca} with $u_0\in B$, $\rho\in(0,1)$ and $\varepsilon\ge0$. Then $\bigl(J_\gamma(u_k)\bigr)$ is nonincreasing and
\begin{equation}
  \sum_k\eta_k\;\le\;\frac{J_\gamma(u_0)}{1-\rho},
  \qquad\text{hence}\qquad
  \min_{k<K}\eta_k\;\le\;\frac{J_\gamma(u_0)}{(1-\rho)K},
  \label{eq:eta-sum}
\end{equation}
where the sum runs over the iterations at which Algorithm~\ref{alg:dca} does not stop. If $\varepsilon>0$, Algorithm~\ref{alg:dca} stops after fewer than $J_\gamma(u_0)/\bigl((1-\rho)\varepsilon\bigr)$ descent steps, at a point $u_k$ with $\theta_k\le\varepsilon$. If $\varepsilon=0$, either some $u_k$ satisfies \eqref{eq:sis_opt_box}, or $(u_k)$ is infinite, has at least one limit point, and every limit point satisfies \eqref{eq:sis_opt_box}.
\end{teo}

\begin{proof}
See Appendix~\ref{app:proofs-method}.
\end{proof}

The proof of Proposition~\ref{prop:descent} uses no strong convexity of $G_\gamma$ or $H_\gamma$, and the decrease required for Theorem~\ref{teo:convergence} is provided by \eqref{eq:model-decrease}. Without strong convexity, our analysis does not guarantee that $\|u_{k+1}-u_k\|$ tends to zero, and for this reason Algorithm~\ref{alg:dca} stops on $\eta_k$ and not on the length of the step. The proof uses only the majorization of $J_\gamma$ by the convex model and the fact that the two certificates in \eqref{eq:relstop} share the dual value. It therefore applies without change to any other scheme whose subproblem majorizes $J_\gamma$ and is solved by Algorithm~\ref{alg:pdhgm}, e.g., to IRL1 (Appendix~\ref{app:competitors}).

In practice, Algorithm~\ref{alg:pdhgm} is also capped at a maximal number of iterations. A capped subproblem only gives $J_\gamma(u_{k+1})\le J_\gamma(u_k)-\eta_k+\mathcal{G}_k(\hat u_{k+1},r_{k+1},s_{k+1})$, so $J_\gamma$ may increase. If only finitely many are capped, Theorem~\ref{teo:convergence} holds from the last capped step on, with $J_\gamma(u_0)$ replaced by the energy at that step.

\section{Numerical experiments}
\label{sec:experiments}
In this section, we evaluate Algorithm~\ref{alg:dca} on denoising, deblurring, inpainting and tomography.
\subsection{Setup}
\label{sec:setup}
\label{sec:protocol}
We tested the method on the thirteen problems of Table~\ref{tab:setup}, which combine four objects with the operators of Section~\ref{sec:introduction}. Two objects are piecewise constant and lie in the model class of \eqref{eq:tvq_problem}: the modified \textsc{Shepp--Logan} phantom \cite{shepp1974fourier} and a \textsc{cartoon} phantom of our own, with bars, discs and steps down to the size of a pixel and to the noise level. The other two are the photographs \textsc{kodim23} \cite{kodak} and butterfly \cite{vanzalinge2019butterfly}, which contain texture and smooth intensity variations. All objects are grayscale images with intensities in $[0,1]$, so they belong to $B$. As in \cite{hintermuller2013nonconvex}, tomography is applied to the \textsc{Shepp--Logan} phantom only, at $64\times64$ pixels. The data is $f=A u^\dagger+e$, with Gaussian noise $e$ of standard deviation $\sigma_{\rm noise}$ and a fixed seed per problem, generated with the same discrete operator used for the reconstruction (Figure~\ref{fig:data}). The construction of the operators and the preparation of the images are described in Appendix~\ref{app:repro}.

Unless stated otherwise, we used $q=0.5$ and $\gamma=100$, so that $\delta_\varphi=7.07$ as in Figure~\ref{fig:huber_smoothing}. For total variation, for the concave model and for the model of Hintermüller and Wu, we selected $\beta$ separately for every object, operator and exponent by the largest PSNR on a geometric grid of ratio $\sqrt2$, whose range is given in Table~\ref{tab:setup}; no selected $\beta$ lies at an endpoint of its grid. Every nonconvex method starts from the total variation solution for the same data and the same $\beta$, computed by Algorithm~\ref{alg:pdhgm} with zero tilt. Algorithm~\ref{alg:pdhgm} uses the rule \eqref{eq:relstop} with $\rho=0.5$ and a cap of $40000$ iterations per subproblem, with the step sizes of Algorithm~\ref{alg:dca} and the balance \eqref{eq:balance}. The first balance is the estimate behind \eqref{eq:balance} on the total variation problem, which is solved with the balance $1$ from a point $u_{\rm init}$ and zero dual variables (Appendix~\ref{app:repro}) and returns $(u_0,r,s)$: we set $c_0=\|u_0-u_{\rm init}\|/\|(r,\delta_\varphi s)\|$, kept in $[10^{-2},10^{2}]$, where $s$ is scaled by $\delta_\varphi$ since the subproblems carry the weight $\lambda=\beta\delta_\varphi$ and the total variation problem the weight $\beta$. Algorithm~\ref{alg:dca} stops once $\eta_k\le\varepsilon$, with $\varepsilon=10^{-4}J_\gamma(u_0)$ for the selection of $\beta$ and $10^{-5}J_\gamma(u_0)$ for the reported reconstructions. The other methods use the stopping rules of their authors, and our implementation of the method of Hintermüller and Wu was validated on two problems of \cite{hintermuller2013nonconvex} (Appendix~\ref{app:competitors}). No subproblem of the reported runs reached the cap.

We measured the quality of a reconstruction $u$ by the PSNR and by the SSIM with the parameters of \cite{wang2004image}. Both are averages over the image and do not separate a loss of contrast at the edges, which a concave penalty is meant to reduce, from other errors. We therefore also measured the recovered edge contrast
\begin{equation}
  \chi(u) := \operatorname{median}_{i\in S} \frac{\phi_i(u)}{\phi_i(u^\dagger)},
  \qquad
  S := \bigl\{ i : \phi_i(u^\dagger)\ge Q_{0.99}(\phi(u^\dagger)) \bigr\},
  \label{eq:contrast}
\end{equation}
i.e., the median recovered jump height on the strongest one percent of the true gradients, where $\chi=1$ means exact recovery and $\chi<1$ means shrinkage. Here, $Q_{0.99}$ denotes the $0.99$-quantile. Following \cite[Table~5]{hintermuller2013nonconvex}, we also report the fraction of pixels with $\phi_i(u)>10^{-3}$, which measures the sparsity of the gradient, and the unsmoothed energy $J(u)$ of Lemma~\ref{lem:smoothing-error}. The implementation is described in Appendix~\ref{app:repro}.

\begin{table}[h]
  \centering
  \caption{Test problems: objects, operators, noise level $\sigma_{\rm noise}$ and the geometric grid of ratio $\sqrt2$ on which $\beta$ is selected. Every object is combined with the three image-domain operators, and \textsc{Shepp--Logan} also with tomography.}
  \label{tab:setup}
  \begin{tabular}{llll}
    \toprule
    object & class & size & source \\
    \midrule
    \textsc{Shepp--Logan} & piecewise constant & $256\times256$, $64\times64$\,$^\dagger$ & \cite{shepp1974fourier} \\
    \textsc{cartoon} & piecewise constant & $256\times256$ & this work \\
    parrots (kodim23) & photograph & $256\times256$ & \cite{kodak} \\
    butterfly & photograph & $256\times256$ & \cite{vanzalinge2019butterfly} \\
    \midrule
    operator & parameters & $\sigma_{\rm noise}$ & $\beta$ grid \\
    \midrule
    denoising & $A=I$ & $0.05$ & $[10^{-3},\,0.256]$ \\
    deblurring & Gaussian $7\times7$, std.\ $1.5$\,$^\ddagger$ & $0.01$ & $[2.5\cdot10^{-4},\,0.064]$ \\
    inpainting & $50\%$ of the pixels\,$^\S$ & $0.05$ & $[10^{-3},\,0.256]$ \\
    tomography & $15$ angles, step $12^\circ$ & $0.05$\,$^\ast$ & $[2.5\cdot10^{-6},\,6.4\cdot10^{-4}]$ \\
    \bottomrule
  \end{tabular}

  \vspace{2pt}
  \noindent\begin{minipage}{\linewidth}
    \raggedright\footnotesize $^\dagger$ For tomography, as in
    \cite[Table~6]{hintermuller2013nonconvex}. $^\ddagger$ Periodic boundary conditions. $^\S$ Observed pixels drawn uniformly at random. $^\ast$ Parallel beam; in the units of line integrals over pixels of unit size, before $A$ is divided by its norm bound.
  \end{minipage}
\end{table}

\begin{figure}[h]
  \centering
  \includegraphics{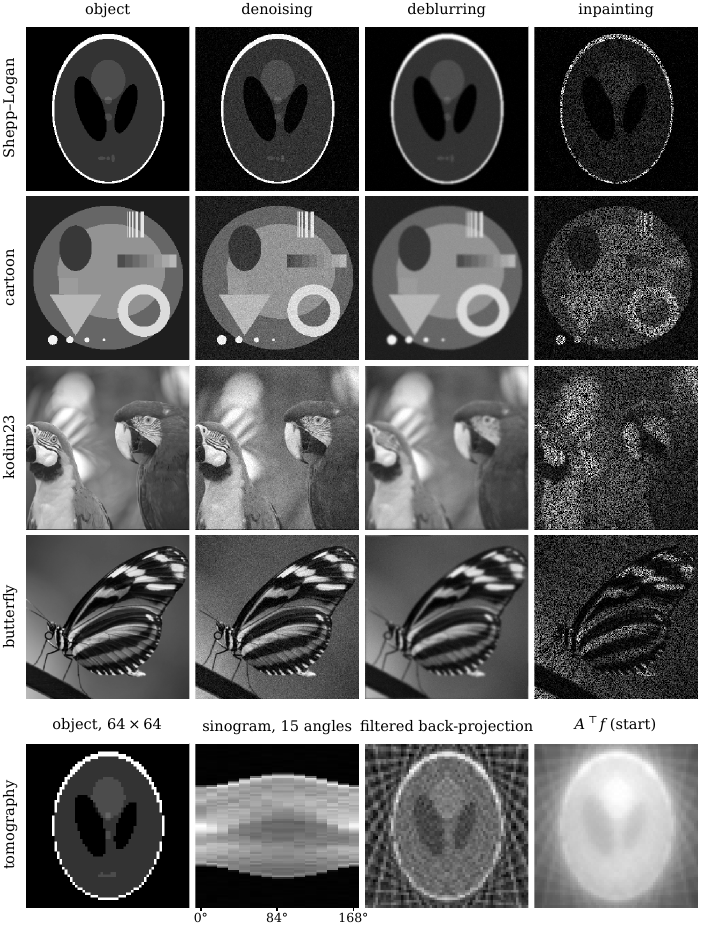}
  \caption{The test problems of Table~\ref{tab:setup}. Top: one row per object, with the object and the data $f$ for denoising, deblurring and inpainting on $[0,1]$; unobserved pixels are black. Bottom: tomography of the $64\times64$ \textsc{Shepp--Logan} phantom, with the object, the sinogram with one column per angle, its filtered back-projection on its own range, and the back-projection $A^\top f$ scaled to its maximum, whose projection onto $B$ starts every run.}
  \label{fig:data}
\end{figure}

\subsection{How the splitting acts on an image}
\label{sec:mechanism}
\begin{figure}[!tbp]
  \centering
  \includegraphics[width=\textwidth]{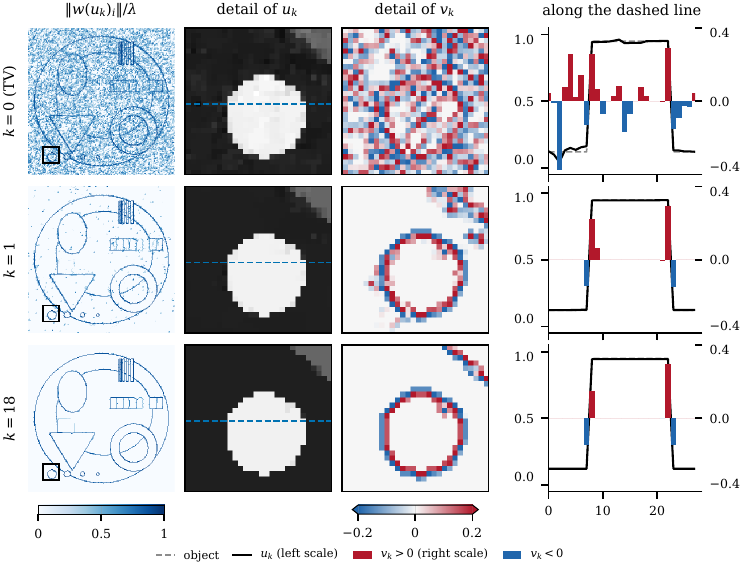}
  \caption{The splitting on an image: denoising of the cartoon phantom with $q=0.5$, $\gamma=100$, $\beta=0.0226$ and $\lambda=0.16$. Rows: the total variation solution for the same $\beta$ ($k=0$), the first outer step, and the point returned by Algorithm~\ref{alg:dca} ($k=18$). Columns: the normalized magnitude $\|w(u_k)_i\|/\lambda$, with the $28\times28$ region of the other columns marked by the black square; the iterate $u_k$ and the tilt $v_k$ in that region; and, along the dashed line, the object and $u_k$ (left scale) and the tilt (bars, right scale).}
  \label{fig:mechanism}
\end{figure}

We first illustrate Lemma~\ref{lem:H_dif} and \eqref{eq:tilt-bound} on the denoising of the cartoon phantom with $q=0.5$ and $\gamma=100$. The rule of Section~\ref{sec:protocol} selected $\beta=0.0226$, so that every subproblem has the weight $\lambda=\beta\delta_\varphi=0.16$, and Algorithm~\ref{alg:dca} stopped after $18$ outer steps. Figure~\ref{fig:mechanism} shows the vector field \eqref{w}, the iterate $u_k$ and the tilt $v_k=\mathbb{K}^\top w(u_k)$ at the total variation solution ($k=0$), after the first outer step and at the returned point.

The first column shows the normalized magnitude $\|w(u_k)_i\|/\lambda$. At every pixel whose gradient magnitude does not exceed the knot $1/\gamma$ we have $w(u_k)_i=0$, and at these pixels the subproblem is the total variation problem with weight $\lambda$. At the total variation solution, the remaining noise places $58\%$ of the gradient magnitudes above the knot, so $w(u_0)$ is spread over the whole image. This fraction is $8.6\%$ after the first outer step and $4.9\%$ at the returned point, where only the edges of the object remain. The bound $\|w(u_k)_i\|\le\lambda$ is almost attained at the strongest edges, where the ratio is $0.936$.

Since $w(u_k)_i$ points in the direction of the gradient $(\mathbb{K}u_k)_i$, the tilt is positive on the bright side of each edge and negative on its dark side, as the profiles of the fourth column show at each crossing of the disc. The linear term $-\scalar{v_k}{u}$ in \eqref{eq:subproblem} therefore favors a larger jump. If the gradient at pixel $i$ keeps its direction, the total variation term and the linear term add up to $(\lambda-\|w(u_k)_i\|)\,\phi_i(u)$, so at the strongest edges the subproblem charges only about $6.5\%$ of the total variation penalty for the height of the jump. At the total variation solution, the noise makes the tilt change sign all over the detail. After one outer step, the tilt is confined to the edge and the noise is gone from the disc and its surroundings. Along the run, the PSNR increases from $33.58$ to $43.89$~dB and the edge contrast \eqref{eq:contrast} from $0.915$ to $0.991$.

\subsection{Certificate and stationarity}
\label{sec:gap_stationarity}
In this subsection, we test Propositions~\ref{prop:gap} and~\ref{prop:descent} on the \textsc{Shepp--Logan} phantom, for the four operators and the values of $\beta$ selected for Algorithm~\ref{alg:dca}.

\begin{figure}[tbp]
  \centering
  \includegraphics[width=\textwidth]{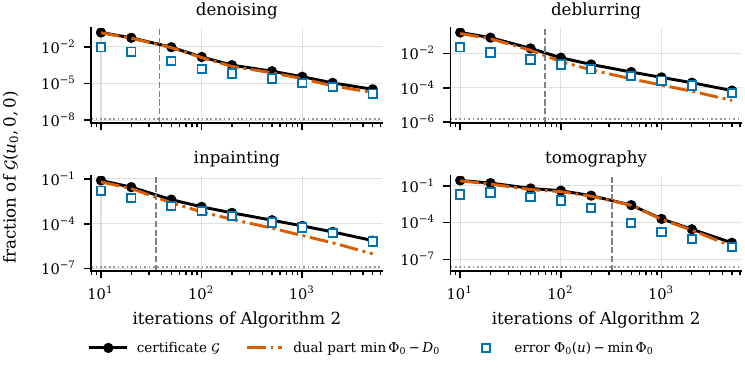}
  \caption{Duality certificate of Proposition~\ref{prop:gap} for the first subproblem on the \textsc{Shepp--Logan} phantom, one panel per operator, against the iterations of Algorithm~\ref{alg:pdhgm} and relative to the certificate at the start. Filled markers: the certificate $\mathcal G_0(u^l,r^l,s^l)$. Hollow markers: the error $\Phi_0(u^l)-\min\Phi_0$ it bounds. Dash-dotted: the distance $\min\Phi_0-D_0(r^l,s^l)$ of the dual value from the minimum. Both are estimated from a reference solve whose certified gap is the dotted line. The dashed vertical line marks the stop of the rule \eqref{eq:relstop} with $\rho=0.5$.}
  \label{fig:certificate}
\end{figure}

For Proposition~\ref{prop:gap}, we solved the first subproblem, with the tilt $v_0$ at the total variation solution, by Algorithm~\ref{alg:pdhgm} from $u_0$ and zero dual variables. The certificate splits as
\[
  \mathcal G_0(u,r,s)=\bigl(\Phi_0(u)-\min\Phi_0\bigr)+\bigl(\min\Phi_0-D_0(r,s)\bigr),
\]
the error it bounds plus the distance of the dual value from the minimum. We estimated both parts with a reference solve of $60\,000$ iterations, and Figure~\ref{fig:certificate} shows them relative to the certificate at the start. The certificate lies above the error at every iteration and for every operator. The difference is largest during the first iterations, more than one order of magnitude for denoising and tomography, where the certificate consists almost entirely of the dual part because the dual variables start at zero. For deblurring and inpainting, the error approaches the certificate as the iterations proceed. The rule \eqref{eq:relstop} with $\rho=0.5$ stops this subproblem after $38$, $70$, $35$ and $320$ iterations for denoising, deblurring, inpainting and tomography.

\begin{figure}[tbp]
  \centering
  \includegraphics[width=\textwidth]{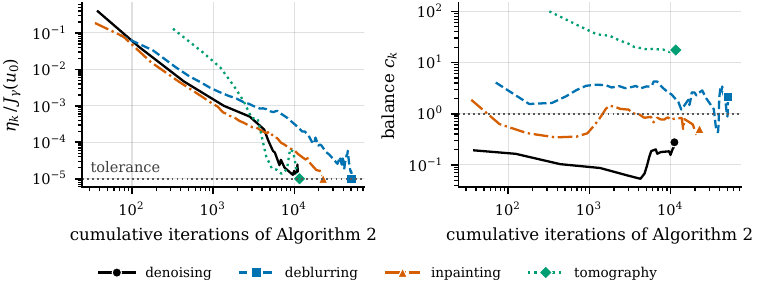}
  \caption{Algorithm~\ref{alg:dca} on the four problems of Figure~\ref{fig:certificate}, under the rule~\eqref{eq:relstop} with $\rho=0.5$, against the cumulative iterations of Algorithm~\ref{alg:pdhgm}. Left: the certified measure $\eta_k$ relative to $J_\gamma(u_0)$, with the tolerance $10^{-5}$ dotted. Right: the balance $c_k$ of \eqref{eq:balance} used in each subproblem, with $c=1$ dotted. Markers: the stop of each run.}
  \label{fig:stationarity}
\end{figure}

For Proposition~\ref{prop:descent}, all four runs of Algorithm~\ref{alg:dca} reached $\eta_k\le10^{-5}J_\gamma(u_0)$, after $1.1\cdot10^4$ to $5.1\cdot10^4$ iterations of Algorithm~\ref{alg:pdhgm} (Figure~\ref{fig:stationarity}, left). The measure $\eta_k$ increases at some steps of every run, which Proposition~\ref{prop:descent} does not exclude. At every outer step, the decrease of $J_\gamma$ was at least $1.11$ times the guaranteed $(1-\rho)\eta_k$ (not shown).

The balance \eqref{eq:balance} (Figure~\ref{fig:stationarity}, right) differs between the operators by more than two orders of magnitude, from about $0.2$ for denoising to $100$ for tomography at the start. It also changes along each run, for tomography from $100$ to $18$.

\subsection{The smoothing parameter $\gamma$}
\label{sec:gamma}
\begin{figure}[tbp]
  \centering
  \includegraphics[width=\textwidth]{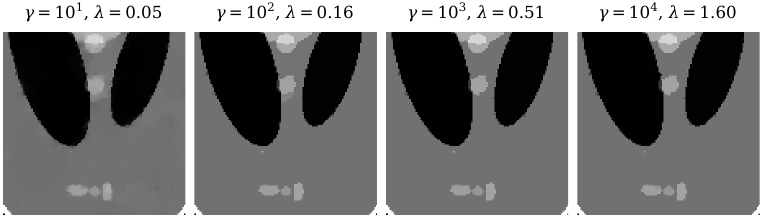}
  \caption{The smoothing parameter on the denoising of the Shepp--Logan phantom, with
  $q=0.5$ and $\beta=0.0226$ fixed: the central region of the reconstruction from the total
  variation solution for four values of $\gamma$, displayed on $[0,0.45]$.}
  \label{fig:gamma}
\end{figure}

We varied $\gamma$ between $10$ and $10^4$ on the denoising of the \textsc{Shepp--Logan} phantom with $q=0.5$ and $\beta=0.0226$ fixed, so that problem \eqref{eq:tvq_problem} is the same in all runs and only its surrogate changes. Algorithm~\ref{alg:dca} was run from the total variation solution (cold start) and by continuation, i.e., from the point returned for the previous value of $\gamma$. Table~\ref{tab:gamma} reports the results, and Figure~\ref{fig:gamma} shows four of the reconstructions.

\begin{table}[tbp]
  \centering
  \caption{The smoothing parameter: denoising of the Shepp--Logan phantom with $q=0.5$ and $\beta=0.0226$. The error per pixel is $(J-J_\gamma)(u_\gamma)/(\beta m)$, with the bound \eqref{eq:smoothing-error}, and iter.\ counts the iterations of Algorithm~\ref{alg:pdhgm} in the run at each $\gamma$.}
  \label{tab:gamma}
  \small\setlength{\tabcolsep}{3.4pt}
\begin{tabular}{rrrrrrrrrr}
\toprule
 & & \multicolumn{2}{c}{error per pixel} & \multicolumn{3}{c}{cold start} & \multicolumn{3}{c}{continuation} \\
\cmidrule(lr){3-4}\cmidrule(lr){5-7}\cmidrule(lr){8-10}
$\gamma$ & $\lambda$ & bound & measured & $J(u_\gamma)$ & PSNR & iter. & $J(u_\gamma)$ & PSNR & iter. \\
\midrule
$10^{1}$ & 0.051 & $1.12\cdot10^{-1}$ & $1.22\cdot10^{-2}$ & 133.83 & 43.62 & 2780 & 133.83 & 43.62 & 2780 \\
$3.16\cdot10^{1}$ & 0.090 & $6.29\cdot10^{-2}$ & $1.55\cdot10^{-3}$ & 120.58 & 46.01 & 5200 & 120.61 & 45.96 & 3256 \\
$10^{2}$ & 0.160 & $3.54\cdot10^{-2}$ & $4.23\cdot10^{-4}$ & 119.55 & 45.81 & 11119 & 119.33 & 46.24 & 9471 \\
$3.16\cdot10^{2}$ & 0.285 & $1.99\cdot10^{-2}$ & $1.15\cdot10^{-4}$ & 119.33 & 45.70 & 11501 & 119.20 & 46.28 & 14480 \\
$10^{3}$ & 0.506 & $1.12\cdot10^{-2}$ & $5.15\cdot10^{-5}$ & 119.32 & 45.67 & 25520 & 119.19 & 46.29 & 8269 \\
$3.16\cdot10^{3}$ & 0.900 & $6.29\cdot10^{-3}$ & $2.07\cdot10^{-5}$ & 119.30 & 45.68 & 64509 & 119.15 & 46.29 & 8782 \\
$10^{4}$ & 1.600 & $3.54\cdot10^{-3}$ & $1.54\cdot10^{-5}$ & 119.35 & 45.55 & 95865 & 119.15 & 46.29 & 5974 \\
\bottomrule
\end{tabular}
\end{table}

Over this range, the bound \eqref{eq:smoothing-error} on the error per pixel decreases by a factor of $32$, and the subproblem weight $\lambda$ increases by the same factor. The measured error stays below the bound. Beyond $\gamma=100$, the unsmoothed energy $J(u_\gamma)$ changes by at most $0.25$, i.e., by $0.2\%$, and the reconstructions in Figure~\ref{fig:gamma} cannot be told apart. At $\gamma=10$, the edges of the small ellipses are blurred and some noise remains in the background.

The cost, on the other hand, grows with $\gamma$. From $\gamma=100$ to $10^4$, the iterations of Algorithm~\ref{alg:pdhgm} grow from $1.1\cdot10^4$ to $9.6\cdot10^4$, while $J(u_\gamma)$ decreases only by $0.2$. Continuation reaches a lower $J$ for every $\gamma\ge10^2$, and its seven runs up to $\gamma=10^4$ take $45\%$ fewer iterations together than the cold start at that value. A moderate $\gamma$ therefore recovers almost all of the decrease of $J$ at a fraction of the cost of a large one.

\subsection{The exponent $q$}
\label{sec:exponent}
\begin{figure}[tbp]
  \centering
  \includegraphics[width=\textwidth]{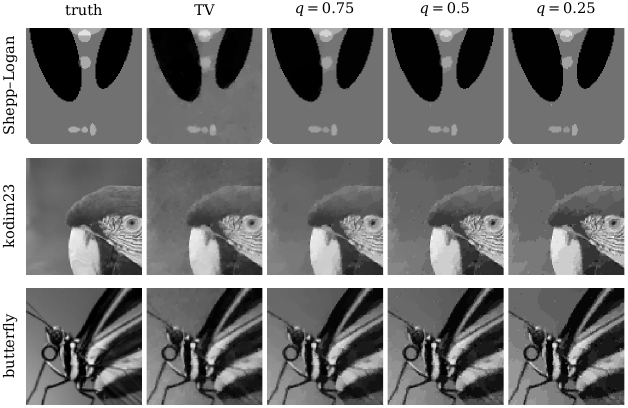}
  \caption{The exponent on the denoising of the Shepp--Logan phantom (top, a region of $112\times112$ pixels displayed on $[0,0.45]$), kodim23 (middle, $96\times96$ pixels) and butterfly (bottom, $96\times96$ pixels), with $\gamma=100$ and $\beta$ selected for every value of $q$.}
  \label{fig:exponent}
\end{figure}

Table~\ref{tab:exponent} and Figure~\ref{fig:exponent} compare total variation with the exponents $q\in\{0.75,0.5,0.25\}$, as in \cite[Table~5]{hintermuller2013nonconvex}, on the denoising of the \textsc{Shepp--Logan} phantom and of the two photographs, with $\gamma=100$ and $\beta$ selected for each value of $q$.

On the phantom, the PSNR, the SSIM, the edge contrast and the sparsity improve as $q$ decreases. The PSNR at $q=0.5$ and $0.25$ is $6.7$ and $7.2$~dB above that of total variation, and the fraction of pixels with a gradient above $10^{-3}$ decreases from $0.264$ to $0.034$, against $0.033$ for the phantom itself.

On the photographs, the edge contrast also increases, from $0.80$ to $0.96$ on kodim23 and from $0.91$ to $1.01$ on butterfly. At $q=0.75$, the PSNR is within $0.2$~dB of that of total variation and the SSIM is not lower. At $q=0.5$ and $0.25$, the PSNR decreases by up to $2.0$~dB and the SSIM by up to $0.07$. The crops in Figure~\ref{fig:exponent} indicate the reason. A small $q$ accentuates the staircase effect of total variation, replacing smooth intensity variations by piecewise constant patches, and keeps isolated noisy pixels as dark dots. The sparsity shows the same effect, since at $q=0.25$ only $0.23$ and $0.25$ of the pixels have a gradient above $10^{-3}$, against $0.99$ and $0.95$ for the photographs themselves.

In every case, $J_\gamma$ is lower at the solution computed for that value of $q$ than at the total variation solution, as observed in \cite[Table~5]{hintermuller2013nonconvex}.

\begin{table}[tbp]
  \centering
  \caption{The exponent, denoising with $\gamma=100$ and $\beta$ selected for every value of $q$. $\chi$ is the edge contrast \eqref{eq:contrast}, sparsity the fraction of pixels with $\phi_i(u)>10^{-3}$ ($0.033$, $0.990$ and $0.952$ for the three objects themselves), and the last two columns the energy of each value of $q$ at its own solution and at the total variation solution.}
  \label{tab:exponent}
  \small\setlength{\tabcolsep}{5pt}
\begin{tabular}{llrrrrrrr}
\toprule
object & $q$ & $\beta$ & PSNR & SSIM & $\chi$ & sparsity & $J_\gamma(\hat u_q)$ & $J_\gamma(\hat u_{\rm TV})$ \\
\midrule
Shepp--Logan & TV & $4.53\cdot10^{-2}$ & 39.12 & 0.962 & 0.933 & 0.264 & \textemdash & \textemdash \\
 & 0.75 & $3.20\cdot10^{-2}$ & 43.68 & 0.985 & 0.984 & 0.056 & 129.74 & 132.85 \\
 & 0.5 & $2.26\cdot10^{-2}$ & 45.81 & 0.990 & 0.995 & 0.037 & 118.92 & 130.64 \\
 & 0.25 & $2.26\cdot10^{-2}$ & 46.34 & 0.992 & 0.998 & 0.034 & 123.86 & 167.89 \\
\midrule
kodim23 & TV & $3.20\cdot10^{-2}$ & 32.68 & 0.886 & 0.803 & 0.763 & \textemdash & \textemdash \\
 & 0.75 & $2.26\cdot10^{-2}$ & 32.52 & 0.887 & 0.871 & 0.501 & 131.42 & 135.87 \\
 & 0.5 & $1.13\cdot10^{-2}$ & 31.68 & 0.846 & 0.946 & 0.353 & 124.47 & 136.31 \\
 & 0.25 & $1.13\cdot10^{-2}$ & 30.64 & 0.816 & 0.960 & 0.233 & 165.65 & 234.06 \\
\midrule
butterfly & TV & $3.20\cdot10^{-2}$ & 32.48 & 0.916 & 0.909 & 0.741 & \textemdash & \textemdash \\
 & 0.75 & $2.26\cdot10^{-2}$ & 32.43 & 0.941 & 0.953 & 0.461 & 146.25 & 151.21 \\
 & 0.5 & $1.13\cdot10^{-2}$ & 31.66 & 0.907 & 0.987 & 0.360 & 132.28 & 144.19 \\
 & 0.25 & $1.13\cdot10^{-2}$ & 30.87 & 0.897 & 1.008 & 0.254 & 173.94 & 238.17 \\
\bottomrule
\end{tabular}
\end{table}

\FloatBarrier
\subsection{Reconstruction across operators and methods}
\label{sec:operators}
We applied both models to the thirteen test problems of Table~\ref{tab:setup}, with $q=0.5$ and $\gamma=100$ for the concave model. All methods for the concave model use the $\beta$ selected for Algorithm~\ref{alg:dca} and start from the total variation solution for that $\beta$. We compared Algorithm~\ref{alg:dca} with three other methods.

IRL1 (own) is the method of \cite{ochs2015iteratively} with the surrogate and the stopping rules of its authors. The surrogate is the shifted power, to which we gave the same slope at the origin as $\varphi_\gamma^q$. Each subproblem stops at the first decrease of the energy, and the run stops once the relative decrease falls below $10^{-7}$. The two surrogates are of a different kind: the shifted power is smooth and modifies $t^q$ for every $t$, whereas $\varphi_\gamma^q$ is linear up to the knot $1/\gamma$.

Algorithm~\ref{alg:dca} (Ochs rules) is Algorithm~\ref{alg:dca} with these two stopping rules in place of its own, so it differs from Algorithm~\ref{alg:dca} only in when it stops. In what follows, we call the rules of Algorithm~\ref{alg:dca} \emph{certified}, since they rest on the duality gap, and those of \cite{ochs2015iteratively} \emph{heuristic}. Finally, the method of Hintermüller and Wu (HW) uses its own quadratic smoothing and has no box constraint (Appendix~\ref{app:competitors}).

\begin{figure}[tbp]
  \centering
  \includegraphics[width=\textwidth]{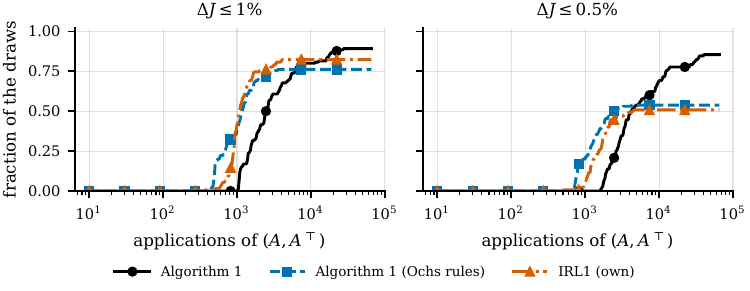}
  \caption{Data profiles of the methods of Table~\ref{tab:comparison-methods} over the $130$ draws of the thirteen test problems: the fraction of the draws on which a method has brought $\Delta J=(J-E^*)/E^*$ below $1\%$ (left) and $0.5\%$ (right) within the work on the horizontal axis, read off its outer iterates.}
  \label{fig:comparison-profile}
\end{figure}

Since each method minimizes a different surrogate, we compared them on the unsmoothed energy $J$, which is common to all. We report $\Delta J=(J-E^*)/E^*$, where $E^*$ is the lowest value found by any method and $J$ is evaluated at $\proj_B(u)$, together with the work in applications of $(A,A^\top)$. Each problem was run on ten draws of the noise and, for inpainting, of the mask, at the $\beta$ selected on the first draw. Table~\ref{tab:comparison-methods} and Figure~\ref{fig:comparison-profile} show on how many of the $130$ draws each method brings $\Delta J$ below $1\%$ and $0.5\%$, and after how much work.

HW is not included in this comparison. Its smoothing leaves small nonzero gradients in the flat regions, and $t^q$ penalizes them heavily. A larger $\Delta J$ would then reflect the smoothing of HW and not the quality of its solver, so we compare HW with Algorithm~\ref{alg:dca} separately, on the smoothed model of each method (Table~\ref{tab:comparison-hw}).

\begin{figure}[tbp]
  \centering
  \includegraphics[width=\textwidth]{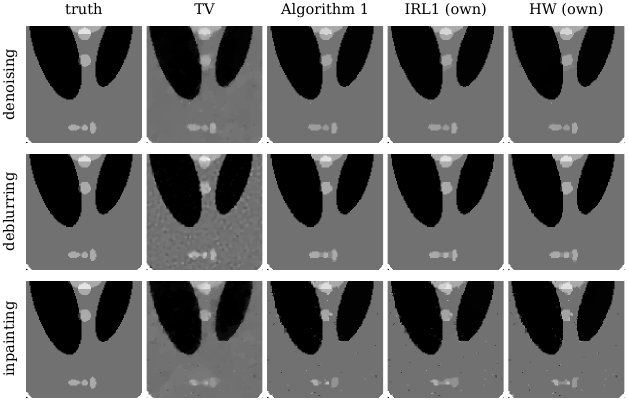}
  \caption{Reconstructions of the Shepp--Logan phantom for the three image-domain
  operators, in the region of Figure~\ref{fig:gamma} displayed on $[0,0.45]$: total
  variation, and the concave model by Algorithm~\ref{alg:dca}, IRL1 (own) and HW (own).
  For tomography, see Figure~\ref{fig:comparison-panels-tomography}.}
  \label{fig:comparison-panels}
\end{figure}

\begin{figure}[tbp]
  \centering
  \includegraphics[width=0.98\textwidth]{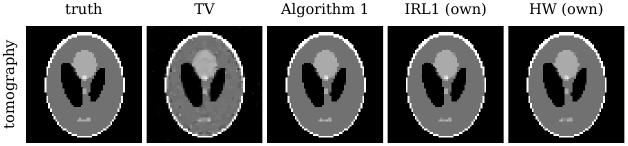}
  \caption{Reconstructions from the tomographic data of the $64\times64$ Shepp--Logan
  phantom, by the methods of Figure~\ref{fig:comparison-panels}: the whole image,
  displayed on $[0,0.45]$ with its pixels unsmoothed.}
  \label{fig:comparison-panels-tomography}
\end{figure}

\paragraph{The two models.} On the piecewise constant objects, the concave model recovers the edges better for every operator. Its PSNR is higher by up to $11.7$~dB for denoising and deblurring and by $14.7$~dB for tomography, and lower by up to $1.2$~dB for inpainting (Table~\ref{tab:comparison-quality} and Figures~\ref{fig:comparison-panels} and~\ref{fig:comparison-panels-tomography}). On the photographs, the edge contrast also increases, while the PSNR decreases by $0.8$ to $2.6$~dB, as in Section~\ref{sec:exponent}.

\paragraph{The stopping rules.} The rules of \cite{ochs2015iteratively} are cheap: the two methods that use them reach the $1\%$ target after $0.5\cdot10^3$ to $1.3\cdot10^3$ applications in median, against $1.5\cdot10^3$ to $7.1\cdot10^3$ for Algorithm~\ref{alg:dca} (Table~\ref{tab:comparison-methods}). Their outer rule, however, measures no distance to a critical point, and the run can stop early. Over the $130$ draws, Algorithm~\ref{alg:dca} (Ochs rules) reaches the $1\%$ target on $99$ and IRL1 (own) on $107$, and both stop up to $3.9\%$ above $E^*$. With its certified rules, Algorithm~\ref{alg:dca} reaches it on $116$ draws, among them every draw of denoising, inpainting and tomography, where it stops at most $0.29\%$ above $E^*$. The $0.5\%$ target separates the rules further, with $111$ draws against $70$ and $66$. Figure~\ref{fig:comparison-profile} shows the same over all draws: the profiles of the heuristic rules rise first, and the profile of Algorithm~\ref{alg:dca} ends highest.

\begin{table}[tbp]
  \centering
  \caption{Methods for the concave model over ten draws per test problem, grouped by operator. $\Delta J=(J-E^*)/E^*$ at the stop, in percent, as median and maximum over the draws; for the targets $1\%$ and $0.5\%$, the number of draws on which the method reaches the target and, in parentheses, the median work of those draws in thousands of applications of $(A,A^\top)$.}
  \label{tab:comparison-methods}
\begin{tabular}{lccc}
\toprule
 & $\Delta J$ at the stop & \multicolumn{2}{c}{draws reaching, and median work to} \\
\cmidrule(lr){3-4}
 & median / max (\%) & $\Delta J\le1\%$ & $\Delta J\le0.5\%$ \\
\midrule
\multicolumn{4}{l}{\textit{denoising}, 40 draws} \\
\quad Algorithm~\ref{alg:dca} & 0.00 / 0.13 & 40/40 (1.5) & 40/40 (2.6) \\
\quad Algorithm~\ref{alg:dca} (Ochs rules) & 0.29 / 3.89 & 36/40 (0.5) & 28/40 (0.8) \\
\quad IRL1 (own) & 0.79 / 3.75 & 25/40 (1.1) & 8/40 (1.8) \\
\midrule
\multicolumn{4}{l}{\textit{deblurring}, 40 draws} \\
\quad Algorithm~\ref{alg:dca} & 0.47 / 1.43 & 26/40 (7.1) & 21/40 (12.2) \\
\quad Algorithm~\ref{alg:dca} (Ochs rules) & 0.96 / 1.52 & 21/40 (1.3) & 13/40 (1.5) \\
\quad IRL1 (own) & 0.00 / 1.50 & 37/40 (0.9) & 37/40 (1.3) \\
\midrule
\multicolumn{4}{l}{\textit{inpainting}, 40 draws} \\
\quad Algorithm~\ref{alg:dca} & 0.05 / 0.26 & 40/40 (2.6) & 40/40 (4.4) \\
\quad Algorithm~\ref{alg:dca} (Ochs rules) & 0.28 / 2.42 & 32/40 (0.7) & 23/40 (1.2) \\
\quad IRL1 (own) & 0.56 / 1.35 & 35/40 (1.2) & 11/40 (2.5) \\
\midrule
\multicolumn{4}{l}{\textit{tomography}, 10 draws} \\
\quad Algorithm~\ref{alg:dca} & 0.19 / 0.29 & 10/10 (2.4) & 10/10 (3.6) \\
\quad Algorithm~\ref{alg:dca} (Ochs rules) & 0.42 / 1.08 & 10/10 (1.0) & 6/10 (1.2) \\
\quad IRL1 (own) & 0.00 / 0.33 & 10/10 (1.0) & 10/10 (1.2) \\
\bottomrule
\end{tabular}
\end{table}

\paragraph{The surrogate.} The draws on which Algorithm~\ref{alg:dca} misses the targets are all of deblurring. There it reaches them on $26$ and $21$ of the $40$ draws and stops up to $1.43\%$ above $E^*$, while IRL1 (own) reaches both on $37$. Algorithm~\ref{alg:dca} (Ochs rules) stops at a similar level, at most $1.52\%$ above, which indicates that this difference is not due to the stopping rules.

\begin{table}[tbp]
  \centering
  \caption{The two models on the thirteen test problems: total variation and the concave model with $q=0.5$ and $\gamma=100$, computed by Algorithm~\ref{alg:dca}, each at its selected $\beta$. $\chi$ is the edge contrast \eqref{eq:contrast}.}
  \label{tab:comparison-quality}
\begin{tabular}{llrrrrrr}
\toprule
 & & \multicolumn{3}{c}{TV} & \multicolumn{3}{c}{$\mathrm{TV}^q$, Algorithm~\ref{alg:dca}} \\
\cmidrule(lr){3-5}\cmidrule(lr){6-8}
object & operator & PSNR & SSIM & $\chi$ & PSNR & SSIM & $\chi$ \\
\midrule
Shepp--Logan & denoising & 39.12 & 0.962 & 0.933 & 45.81 & 0.990 & 0.995 \\
 & deblurring & 33.57 & 0.949 & 0.929 & 45.25 & 0.996 & 0.999 \\
 & inpainting & 28.19 & 0.916 & 0.852 & 26.99 & 0.954 & 0.989 \\
 & tomography & 33.33 & 0.970 & 0.916 & 48.02 & 0.998 & 0.997 \\
\midrule
cartoon & denoising & 38.00 & 0.980 & 0.877 & 43.89 & 0.996 & 0.991 \\
 & deblurring & 33.20 & 0.957 & 0.714 & 35.54 & 0.988 & 0.997 \\
 & inpainting & 31.69 & 0.965 & 0.728 & 31.65 & 0.973 & 0.902 \\
\midrule
kodim23 & denoising & 32.68 & 0.886 & 0.803 & 31.68 & 0.846 & 0.946 \\
 & deblurring & 29.35 & 0.883 & 0.292 & 28.37 & 0.865 & 0.313 \\
 & inpainting & 28.68 & 0.832 & 0.529 & 27.54 & 0.774 & 0.701 \\
\midrule
butterfly & denoising & 32.48 & 0.916 & 0.909 & 31.66 & 0.907 & 0.987 \\
 & deblurring & 30.06 & 0.933 & 0.810 & 27.46 & 0.918 & 1.034 \\
 & inpainting & 28.29 & 0.869 & 0.761 & 26.75 & 0.836 & 0.960 \\
\bottomrule
\end{tabular}
\end{table}

Table~\ref{tab:comparison-hw} compares Algorithm~\ref{alg:dca} with HW, with each solution evaluated on both smoothed models. On $J_\gamma$, the solution of HW is higher on all thirteen problems, by $0.06\%$ to $1.47\%$. On $J_{\rm HW}$, it is lower on eleven problems, by at most $1.13\%$, and higher on the denoising and the inpainting of kodim23, by at most $0.13\%$. The two solutions are therefore close in either energy. At the $\beta$ selected for each method, the PSNR of HW is within $0.21$~dB of Algorithm~\ref{alg:dca} on the twelve image-domain problems and $1.65$~dB lower for tomography.

The reconstructions differ in the flat regions: at most $0.2\%$ of the pixels of HW have a gradient below $10^{-6}$, against $38.5\%$ to $95.9\%$ for Algorithm~\ref{alg:dca}. Regarding the work, HW needs less than Algorithm~\ref{alg:dca} on four problems, namely tomography, the deblurring of the \textsc{Shepp--Logan} phantom and the denoising of the two piecewise constant objects, and $1.1$ to $4.6$ times more on the other nine, with at least twice as much on the photographs.

\begin{table}[tbp]
  \centering
  \caption{HW against Algorithm~\ref{alg:dca}, as medians over the ten draws. Energies of the solution of HW relative to those of Algorithm~\ref{alg:dca}, at the same $\beta$: the model $J_\gamma$ of Algorithm~\ref{alg:dca}, at $\proj_B(u)$, and the model $J_{\rm HW}$ (Appendix~\ref{app:competitors}). Other columns at the $\beta$ selected for each method: PSNR of HW minus that of Algorithm~\ref{alg:dca}, pixels with gradient at most $10^{-6}$, and work of HW over that of Algorithm~\ref{alg:dca}.}
  \label{tab:comparison-hw}
\begin{tabular}{llrrrrrrr}
\toprule
 & & & \multicolumn{2}{c}{HW relative to Alg.~\ref{alg:dca} (\%)} & PSNR & \multicolumn{2}{c}{gradient $\le10^{-6}$ (\%)} & work \\
\cmidrule(lr){4-5}\cmidrule(lr){7-8}
object & operator & $\beta_{\rm HW}/\beta$ & $J_\gamma$ & $J_{\rm HW}$ & $\Delta$ (dB) & Alg.~\ref{alg:dca} & HW & ratio \\
\midrule
Shepp--Logan & denoising & 1.00 & $+0.16$ & $-0.05$ & $-0.21$ & 95.9 & 0.1 & 0.8 \\
 & deblurring & 2.00 & $+0.14$ & $-0.19$ & $+0.10$ & 92.6 & 0.2 & 0.4 \\
 & inpainting & 2.00 & $+0.06$ & $-0.18$ & $-0.05$ & 92.0 & 0.2 & 2.1 \\
 & tomography & 1.41 & $+1.47$ & $-1.08$ & $-1.65$ & 64.6 & 0.0 & 0.2 \\
\midrule
cartoon & denoising & 1.41 & $+0.22$ & $-0.06$ & $-0.05$ & 94.1 & 0.1 & 0.9 \\
 & deblurring & 1.00 & $+0.32$ & $-1.13$ & $+0.06$ & 48.2 & 0.0 & 1.2 \\
 & inpainting & 1.00 & $+0.06$ & $-0.18$ & $+0.09$ & 94.9 & 0.2 & 1.1 \\
\midrule
kodim23 & denoising & 1.00 & $+0.90$ & $+0.07$ & $-0.01$ & 55.8 & 0.0 & 3.1 \\
 & deblurring & 1.00 & $+1.10$ & $-0.61$ & $-0.04$ & 43.5 & 0.0 & 3.1 \\
 & inpainting & 1.41 & $+0.65$ & $+0.13$ & $-0.01$ & 62.1 & 0.0 & 4.1 \\
\midrule
butterfly & denoising & 1.41 & $+0.82$ & $-0.14$ & $-0.04$ & 51.2 & 0.0 & 2.4 \\
 & deblurring & 0.71 & $+1.00$ & $-0.76$ & $-0.11$ & 38.5 & 0.0 & 2.0 \\
 & inpainting & 1.00 & $+0.59$ & $-0.02$ & $-0.03$ & 58.3 & 0.0 & 4.6 \\
\bottomrule
\end{tabular}
\end{table}

\FloatBarrier
\section{Conclusions and limitations}
\label{sec:conclusions}
In this work, we investigated image reconstruction regularized by the quasinormed total variation $\text{TV}^q$ with a general linear forward operator, and its numerical solution by difference-of-convex programming. Since $t^q$ has unbounded slope at the origin, the objective admits no splitting whose convex part is a total variation functional, and we replaced it by a surrogate of finite slope. The resulting objective is the difference of a total variation functional with a single scalar weight and a continuously differentiable convex function. All subproblems are therefore total variation problems with the same weight and differ only in a linear term. We solved them by a primal--dual method that only applies the forward operator and, with the image restricted to its dynamic range, returns a computable duality gap. Stopping each subproblem once its gap falls below a fixed fraction of the gap at the current outer iterate certifies the descent of every outer step, bounds the number of outer steps and yields critical points of the smoothed problem over the dynamic range.

In the experiments, the concave model recovered the edges of the piecewise constant objects better than total variation for every operator, with a PSNR up to $11.7$~dB higher in the image domain, except for inpainting, and $14.7$~dB higher for tomography. On the photographs, the edge contrast also increased, but the PSNR decreased, since a small exponent replaces smooth intensity variations by piecewise constant patches. With its certified rules, Algorithm~\ref{alg:dca} brought the unsmoothed energy within $1\%$ of the best value found on $116$ of the $130$ draws, while the heuristic rules of the other methods stopped up to $4\%$ above it. The remaining draws are all of deblurring, where it stopped at most $1.4\%$ above the best value. At its own regularization parameter, the method of Hintermüller and Wu reached a PSNR within $0.3$~dB of Algorithm~\ref{alg:dca} in the image domain and $1.65$~dB below for tomography.

Some limitations of this study should be mentioned. The regularization parameter was selected by the PSNR against the ground truth, so the results compare the models at their best parameter and do not address its selection in practice. The convergence result concerns critical points of the smoothed problem over the dynamic range, and different methods, or different stopping rules of the same method, stop at different critical points. The certified rules are more expensive: to reach the $1\%$ target, Algorithm~\ref{alg:dca} needed $1.5$ to $7.1$ thousand applications of the operator and its adjoint in median, against $0.5$ to $1.3$ thousand for the heuristic rules. Learning the model parameters from training data by a bilevel strategy, instead of the grid search used here, is a promising research direction, and the extension to vector-valued images and to higher-order regularizers is a matter for further research.

\begin{appendices}
\renewcommand{\theHsection}{app.\thesection}
\renewcommand{\theHequation}{app.\thesection.\theequation}

\section{Proofs for the DC splitting and the smoothing error}
\label{app:proofs-splitting}
\begin{proof}[Proof of Lemma~\ref{lem:H_dif}]
\emph{The scalar function $j_\varphi$.} By \ref{ass:A2} and Lemma~\ref{lem:varphi_prime_nneg_ninc}, $j_\varphi'=\delta_\varphi-\varphi'$ is continuous on $(0,\infty)$, satisfies $0\le j_\varphi'\le\delta_\varphi$ there, and tends to $0$ as $t\downarrow0$. At the origin, $j_\varphi(0)=0$ and
\begin{equation*}
  \frac{j_\varphi(t)-j_\varphi(0)}{t}=\delta_\varphi-\frac{\varphi(t)}{t}\;\xrightarrow[t\downarrow0]{}\;0
\end{equation*}
by \ref{ass:A3}, so the right derivative of $j_\varphi$ at $0$ exists and equals $0=\lim_{t\downarrow0}j_\varphi'(t)$. Hence $j_\varphi\in C^1([0,\infty))$ with $j_\varphi'(0)=0$ and $0\le j_\varphi'\le\delta_\varphi$ on $[0,\infty)$.

\emph{The function $H_\varphi$.} Write
\begin{equation*}
  H_i(u):=j_\varphi\bigl(\phi_i(u)\bigr),\qquad i=1,\dots,m,
\end{equation*}
so that $H_\varphi=\beta\sum_{i=1}^mH_i$. We show that every $H_i$ is differentiable on $\mathbb{R}^n$, with
\begin{equation}
    \nabla H_i(u)=\frac{1}{\beta}\,\mathbb{K}_i^\top w_\varphi(u)_i .
    \label{eq:Hi-grad}
\end{equation}

\emph{Blocks with $\phi_i(u)>0$.} The norm is differentiable on $\mathbb{R}^2\setminus\{0\}$ with gradient $x/\|x\|$, and $j_\varphi$ is differentiable at $\phi_i(u)>0$. By the chain rule, $H_i$ is differentiable at $u$, with
\begin{equation*}
  \nabla H_i(u)=j_\varphi'\bigl(\phi_i(u)\bigr)\,\mathbb{K}_i^\top\frac{(\mathbb{K}u)_i}{\phi_i(u)}
=\frac{1}{\beta}\,\mathbb{K}_i^\top w_\varphi(u)_i,
\end{equation*}
the second equality being the definition \eqref{w} of $w_\varphi$. In particular, for the penalty \eqref{eq:reg} no separate case is needed at the smoothing knot.

\emph{Blocks with $\phi_i(u)=0$.} Here $(\mathbb{K}u)_i=0$, so $H_i(u+h)-H_i(u)=j_\varphi\bigl(\|\mathbb{K}_ih\|\bigr)-j_\varphi(0)$, with $\|\mathbb{K}_ih\|\le\|\mathbb{K}_i\|\,\|h\|$. Since $j_\varphi'(0)=0$, we have $j_\varphi(s)-j_\varphi(0)=o(s)$ as $s\downarrow0$, and therefore $H_i(u+h)-H_i(u)=o(\|h\|)$ as $h\to0$. Hence $H_i$ is differentiable at $u$ with $\nabla H_i(u)=0$, which is \eqref{eq:Hi-grad} under the convention $w_\varphi(u)_i=0$.

\emph{Conclusion.} Summing the blocks, $H_\varphi$ is differentiable on $\mathbb{R}^n$ with
\begin{equation*}
  \nabla H_\varphi(u)=\sum_{i=1}^m\mathbb{K}_i^\top w_\varphi(u)_i=\mathbb{K}^\top w_\varphi(u).
\end{equation*}
It remains to show that $w_\varphi$ is continuous on $\mathbb{R}^n$. On a block with $\phi_i(u)>0$, this follows from the continuity of $\phi_i$, of $j_\varphi'$ on $(0,\infty)$ and of $x\mapsto x/\|x\|$ away from the origin. On a block with $\phi_i(u)=0$, we have $\|w_\varphi(u')_i\|=\beta\,j_\varphi'\bigl(\phi_i(u')\bigr)\to\beta\,j_\varphi'(0)=0$ as $u'\to u$, by the continuity of $\phi_i$ and of $j_\varphi'$ on $[0,\infty)$. Hence $\nabla H_\varphi$ is continuous, and $H_\varphi$ is continuously differentiable on $\mathbb{R}^n$.
\end{proof}

\begin{proof}[Proof of Theorem~\ref{teo:opt_cond}]
Let $u^*$ be a local minimizer of \eqref{eq:Pphi}, i.e., of $J_\varphi=G_\varphi-H_\varphi$, where
\begin{equation*}
  G_\varphi(u)=\frac{1}{2}\|Au-f\|^2+\lambda\sum_{i=1}^m\|(\mathbb{K}u)_i\|,\qquad \lambda=\beta\delta_\varphi .
\end{equation*}
Since $G_\varphi$ and $H_\varphi$ are convex and $H_\varphi$ is finite everywhere, the necessary optimality condition for a local minimizer of a difference-of-convex program \cite{hor} gives
\begin{equation*}
  \partial H_\varphi(u^*)\subseteq\partial G_\varphi(u^*).
\end{equation*}
By Lemma~\ref{lem:H_dif}, $H_\varphi$ is differentiable at $u^*$ and, since it is a finite convex function, its subdifferential reduces to the gradient~\cite{urru}, i.e., $\partial H_\varphi(u^*)=\{\mathbb{K}^\top w_\varphi(u^*)\}$. The condition above therefore reads $\mathbb{K}^\top w_\varphi(u^*)\in\partial G_\varphi(u^*)$.

It remains to compute $\partial G_\varphi(u^*)$. The two terms of $G_\varphi$ are convex, finite and continuous on $\mathbb{R}^n$, so the sum rule of subdifferential calculus holds without further qualification conditions \cite{urru}, and
\begin{equation*}
  \partial G_\varphi(u^*)=\partial\left(\frac{1}{2}\|A\cdot-f\|^2\right)(u^*)+\partial \left(\lambda\sum_{i=1}^m\|\mathbb{K}_i\cdot\|\right)(u^*),
\end{equation*}
The fidelity term is differentiable, so its subdifferential is $\{A^\top(Au^*-f)\}$. For the total variation term, we apply the sum rule to the finite convex functions $u\mapsto\|\mathbb{K}_iu\|$, together with the chain rule for a linear map, which is exact since $\|\cdot\|$ is finite everywhere \cite{urru}, and obtain
\begin{align*}
    \partial(\|\cdot\|\circ\mathbb{K}_i)(u^*)=\mathbb{K}_i^\top\partial\|\cdot\|\bigl((\mathbb{K}u^*)_i\bigr),
    \qquad i=1,\ldots,m.
\end{align*}
Consequently,
\begin{align*}
    \partial G_\varphi(u^*)=\left\{A^\top(Au^*-f)+\lambda\sum_{i=1}^m\mathbb{K}_i^\top z_i: z_i\in \partial \|\cdot\|\bigl((\mathbb{K}u^*)_i\bigr)\right\},
\end{align*}
where, for each $i=1,\ldots,m$,
\begin{align*}
    \partial \|\cdot\|\bigl((\mathbb{K}u^*)_i\bigr)=\left\{\begin{array}{ll}
    \displaystyle\left\{\frac{(\mathbb{K}u^*)_i}{\|(\mathbb{K}u^*)_i\|}\right\}, & (\mathbb{K}u^*)_i\neq 0, \\[10pt]
    \{v\in \mathbb{R}^2:\|v\|\leq 1\}, & (\mathbb{K}u^*)_i=0.
    \end{array}\right.
\end{align*}
Finally, stacking $z=(z_1,\dots,z_m)\in\mathbb{R}^{2m}$ in the same block order as $w_\varphi$, we have $\sum_{i=1}^m\mathbb{K}_i^\top z_i=\mathbb{K}^\top z$. The inclusion $\mathbb{K}^\top w_\varphi(u^*)\in\partial G_\varphi(u^*)$ is then equivalent to the existence of $z_i\in\partial\|\cdot\|\bigl((\mathbb{K}u^*)_i\bigr)$ for all $i=1,\ldots,m$, with
\begin{align*}
    \mathbb{K}^\top w_\varphi(u^*)-A^\top(Au^*-f)-\lambda\mathbb{K}^\top z=0,
\end{align*}
which is \eqref{eq:sis_opt_low}.
\end{proof}

\begin{proof}[Proof of Lemma~\ref{lem:smoothing-error}]
We begin by establishing a linear upper bound for the function $x\mapsto x^q$. Since $q \in (0,1)$, its second derivative $q(q-1)x^{q-2}$ is strictly negative for all $x > 0$, meaning the function is strictly concave. By the first-order characterization of concavity, the graph of the function is globally bounded above by any of its tangent lines. Evaluating the tangent at an arbitrary point $a > 0$, we directly obtain:
\begin{equation}
x^q \le a^q + qa^{q-1}(x-a). \label{eq:tangent-bound}
\end{equation}

We evaluate this tangent line at the specific point $a = q/\gamma$. For this choice, the slope is $qa^{q-1} = q(q/\gamma)^{q-1} = q^q \gamma^{1-q} = \delta_\varphi$, and also we have 
\begin{equation*} 
a^q - \delta_{\varphi}a = a^q-(q^q\gamma^{1-q})a= a^q(1-q) = \left(\frac{q}{\gamma}\right)^q(1-q) = q^q(1-q)\gamma^{-q}. 
\end{equation*}
Substituting $a$ and $\delta_\varphi$ into \eqref{eq:tangent-bound} and using the above equality, we obtain that for any $x \ge 0$:
\begin{align*}
  x^q\leq \left(\dfrac{q}{\gamma}\right)^q + \delta_\varphi(x - \dfrac{q}{\gamma}) = \delta_\varphi x + q^q(1-q)\gamma^{-q},
\end{align*}
then 
\begin{equation}
  x^q-\delta_\varphi x \le q^q(1-q)\gamma^{-q}.
  \label{eq:tangent-bound-eval}
\end{equation}

We now analyze the approximation error $t^q - \varphi_\gamma^q(t)$ over its two distinct subdomains. In the first interval, $t \in [0, 1/\gamma]$, the smoothed function is linear: $\varphi_\gamma^q(t) = \delta_\varphi t$. Using the derived tangent bound \eqref{eq:tangent-bound-eval} evaluated at $x=t$, we rearrange the terms to immediately obtain the upper bound:
\begin{equation*} 
t^q - \delta_\varphi t \le q^q(1-q)\gamma^{-q}. 
\end{equation*}
Note that equality is achieved exactly at the tangency point $t = q/\gamma$. To establish the lower bound in this interval, we must show that $t^q - \delta_\varphi t \ge 0$. Since $q \in (0,1)$, the exponent $(q-1)$ is negative, making the function $t^{q-1}$ strictly decreasing. Thus, for $t \in (0, 1/\gamma]$, it attains its minimum at the right endpoint:
\begin{equation*}
t^{q-1} \ge \left(\frac{1}{\gamma}\right)^{q-1} = \gamma^{1-q}. 
\end{equation*}
Because $q^q < 1$, it follows that $\gamma^{1-q} > q^q\gamma^{1-q} = \delta_\varphi$. Therefore, $t^{q-1} \ge \delta_\varphi$, which upon multiplying by $t$ yields $t^q \ge \delta_\varphi t$, guaranteeing a non-negative error.

For the second interval, $t > 1/\gamma$, the smoothed function is a translation: $\varphi_\gamma^q(t) = x^q$, where we define $x := t - \frac{1-q}{\gamma}$. Notice that $x$ lies strictly to the right of our original tangency point $a$:
$$ x = t - \frac{1-q}{\gamma} > \frac{1}{\gamma} - \frac{1-q}{\gamma} = \frac{q}{\gamma} = a. $$
The error in this region is $t^q - x^q$. Applying the general tangent bound \eqref{eq:tangent-bound} again, but this time centered at $x$ and evaluated at $t$, we have $t^q \le x^q + qx^{q-1}(t-x)$. Substituting the distance $(t-x) = \frac{1-q}{\gamma}$ yields:
\begin{equation*}
t^q - x^q \le qx^{q-1}\left(\frac{1-q}{\gamma}\right). 
\end{equation*}
Since $x > q/\gamma$ and the exponent $(q-1)$ is negative, the decreasing nature of the power function implies $x^{q-1} < (q/\gamma)^{q-1}$. Substituting this into our inequality gives:
\begin{equation*}
t^q - x^q \le q\left(\frac{q}{\gamma}\right)^{q-1} \frac{1-q}{\gamma} = \left(q^q \gamma^{1-q}\right) \frac{1-q}{\gamma} = q^q(1-q)\gamma^{-q}. 
\end{equation*}
The lower bound $t^q - x^q \ge 0$ is trivial since $t > x$. Consequently, for all $t \ge 0$, the point-wise error satisfies $0 \le t^q - \varphi_\gamma^q(t) \le q^q(1-q)\gamma^{-q}$.

To extend this result to the global objective function, we note that the difference between the unsmoothed objective $J(u)$ and the smoothed objective $J_\gamma(u)$ arises solely from the $m$ pixels in the regularization term. Letting $t = \|(\mathbb{K}u)_i\|$, we have:
\begin{equation*}
J(u) - J_\gamma(u) = \beta \sum_{i=1}^m \left( \|(\mathbb{K}u)_i\|^q - \varphi_\gamma^q\left(\|(\mathbb{K}u)_i\|\right) \right). 
\end{equation*}
Since each term in the sum is bounded between $0$ and $q^q(1-q)\gamma^{-q}$, the total difference is bounded by:
$$ 0 \le J(u) - J_\gamma(u) \le \beta \, m \, q^q(1-q)\gamma^{-q} = E(\gamma). $$

Finally, let $u_\gamma$ be any minimizer of $J_\gamma$. For an arbitrary $u \in \R^n$, we can construct a chain of inequalities using the established bound $E(\gamma)$. We know that $J(u_\gamma) \le J_\gamma(u_\gamma) + E(\gamma)$. Because $u_\gamma$ is the global minimum of $J_\gamma$, it holds that $J_\gamma(u_\gamma) \le J_\gamma(u)$. Substituting this yields:
\begin{equation*}
J(u_\gamma) \le J_\gamma(u) + E(\gamma). 
\end{equation*}
Since the approximation error is non-negative, $J_\gamma(u) \le J(u)$, leading to our final conclusion:
$$ J(u_\gamma) \le J(u) + E(\gamma). $$
This proves that the minimizer of the smoothed problem yields an objective value within $E(\gamma)$ of the true infimum of the original problem.
\end{proof}

\section{Proofs for the DC method}
\label{app:proofs-method}

\begin{proof}[Proof of Proposition~\ref{prop:descent}]
\emph{Majorization.} Put $C_k:=\scalar{v_k}{u_k}-H_\gamma(u_k)$. Since $H_\gamma$ is convex and differentiable with $\nabla H_\gamma(u_k)=v_k$ (Lemma~\ref{lem:H_dif}), we have $H_\gamma(u)\ge H_\gamma(u_k)+\scalar{v_k}{u-u_k}$ for every $u\in\mathbb{R}^n$, and therefore
\begin{equation}
  J_\gamma(u)=G_\gamma(u)-H_\gamma(u)\le G_\gamma(u)-\scalar{v_k}{u}+C_k\le\Phi_k(u)+C_k,
  \qquad u\in\mathbb{R}^n,
  \label{eq:majorization}
\end{equation}
with equality at $u=u_k\in B$.

\emph{Duality bound.} By weak duality, $D_k(r_{k+1},s_{k+1})\le\min\Phi_k$, so $\eta_k=\Phi_k(u_k)-D_k(r_{k+1},s_{k+1})\ge\theta_k$. If $\eta_k>\varepsilon$, the triple returned by Algorithm~\ref{alg:pdhgm} fails the second test in \eqref{eq:relstop} and therefore satisfies the first, which, since both certificates share the dual value, reads $\mathcal{G}_k(\hat u_{k+1},r_{k+1},s_{k+1})\le\rho\,\eta_k$. Hence
\begin{align*}
  \eta_k-\theta_k=\min\Phi_k-D_k(r_{k+1},s_{k+1})
  &\le\Phi_k(\hat u_{k+1})-D_k(r_{k+1},s_{k+1})\\
  &=\mathcal{G}_k(\hat u_{k+1},r_{k+1},s_{k+1})\le\rho\,\eta_k ,
\end{align*}
which is \eqref{eq:eta-theta}. Moreover, $u_{k+1}=\hat u_{k+1}\in B$, since Algorithm~\ref{alg:pdhgm} returns a projection onto $B$.

\emph{Descent.} Since $u_{k+1}=\hat u_{k+1}$, \eqref{eq:majorization}, the equality in~\eqref{eq:majorization} at $u_k$ and~\eqref{eq:model-decrease} with $\mathcal{G}_k(u_k,r_{k+1},s_{k+1})=\eta_k$ give
\begin{align*}
  J_\gamma(u_{k+1})\le\Phi_k(\hat u_{k+1})+C_k
  &=J_\gamma(u_k)-\bigl(\Phi_k(u_k)-\Phi_k(\hat u_{k+1})\bigr)\\
  &\le J_\gamma(u_k)-(1-\rho)\,\eta_k ,
\end{align*}
which is \eqref{eq:descent}.
\end{proof}

\begin{proof}[Proof of Theorem~\ref{teo:convergence}]
By assumption $u_0\in B$, so every iterate lies in $B$ by Proposition~\ref{prop:descent}.

\emph{Summability.} By \eqref{eq:descent}, $\bigl(J_\gamma(u_k)\bigr)$ is nonincreasing, and $J_\gamma\ge0$, since the fidelity term is nonnegative and $\varphi_\gamma^q\ge0$. Summing \eqref{eq:descent} over the iterations $k=0,\dots,K-1$ at which Algorithm~\ref{alg:dca} does not stop gives
\begin{equation*}
  (1-\rho)\sum_{k=0}^{K-1}\eta_k\;\le\;J_\gamma(u_0)-J_\gamma(u_K)\;\le\;J_\gamma(u_0),
\end{equation*}
which is \eqref{eq:eta-sum}; in particular, $\eta_k\to0$, and $\theta_k\to0$ by \eqref{eq:eta-theta}. If $\varepsilon>0$ and Algorithm~\ref{alg:dca} does not stop at the iterations $k=0,\dots,K-1$, then $\eta_k>\varepsilon$ for each of them, and the last display gives $K<J_\gamma(u_0)/\bigl((1-\rho)\varepsilon\bigr)$. At the iteration where it stops, $\theta_k\le\eta_k\le\varepsilon$ by Proposition~\ref{prop:descent}.

\emph{Stationarity.} Let $\varepsilon=0$. If Algorithm~\ref{alg:dca} stops at $u_k$, then $\eta_k=0$, so $\theta_k=0$, and $u_k$ satisfies \eqref{eq:sis_opt_box}. If Algorithm~\ref{alg:pdhgm} does not terminate at iteration $k$, then $\theta_k=0$ by Lemma~\ref{lem:inner-stop}, and the same conclusion holds. Otherwise, $(u_k)$ is an infinite sequence in the compact set $B$ and has at least one limit point. Let $u_{k_j}\to\bar u\in B$ and fix $y\in B$. Since $v_k=\nabla H_\gamma(u_k)$ and $\min\Phi_k\le\Phi_k(y)$,
\begin{equation*}
  G_\gamma(u_k)-G_\gamma(y)+\scalar{\nabla H_\gamma(u_k)}{y-u_k}
  =\Phi_k(u_k)-\Phi_k(y)\le\theta_k .
\end{equation*}
Since $G_\gamma$ and $\nabla H_\gamma$ are continuous (Lemma~\ref{lem:H_dif}), letting $k=k_j\to\infty$ and using $\theta_k\to0$ gives $G_\gamma(\bar u)-G_\gamma(y)+\scalar{\nabla H_\gamma(\bar u)}{y-\bar u}\le0$. Since $y\in B$ is arbitrary, $\bar u$ minimizes $G_\gamma+\iota_B-\scalar{\nabla H_\gamma(\bar u)}{\cdot}$, and the argument of Corollary~\ref{cor:opt_box} shows that $\bar u$ satisfies \eqref{eq:sis_opt_box}.
\end{proof}

\section{Implementation of the comparison methods}
\label{app:competitors}
\paragraph{IRL1.} We implemented Algorithm~3 of \cite{ochs2015iteratively} on the
constraint set $B$. At $u_k$, it replaces the penalty $\varphi$ in use by its tangent at $\phi_i(u_k)$ in every block and computes
\begin{equation}
  u_{k+1}\in\argmin_{u}\;\Psi_k(u):=\tfrac12\|Au-f\|^2+\sum_{i=1}^m\omega_i^k\,\phi_i(u)+\iota_B(u),
  \quad \omega_i^k:=\beta\,\varphi'\bigl(\phi_i(u_k)\bigr),
  \label{eq:irl1}
\end{equation}
a total variation problem with pixel-dependent weights. Since $\varphi$ is concave, $J_\varphi(u)-J_\varphi(u_k)\le\Psi_k(u)-\Psi_k(u_k)$ for every $u\in B$. For the comparison with each method's own surrogate, IRL1 is applied to $\varphi(t)=(t+\varsigma)^q-\varsigma^q$, which is admissible with $\delta_\varphi=q\varsigma^{q-1}$.

Problem \eqref{eq:irl1} has the saddle-point form \eqref{eq:saddle_point} with $v_k=0$ and one ball of radius $\omega_i^k$ per pixel in place of the ball of radius $\lambda$ in \eqref{eq:dualball}. We solved it by Algorithm~\ref{alg:pdhgm} with the same step sizes as \eqref{eq:subproblem}, the balance \eqref{eq:balance} and the same $c_0$.

We used two inner stopping rules. The first is the relative rule \eqref{eq:relstop} with $\rho=0.5$ and $\Psi_k$ in place of $\Phi_k$. With this rule, Proposition~\ref{prop:gap}, Lemma~\ref{lem:inner-stop}, Proposition~\ref{prop:descent} and \eqref{eq:eta-sum} hold for IRL1 with $\Psi_k$ and the radii $\omega_i^k$ in place of $\Phi_k$ and $\lambda$, since their proofs use only the separability of $\tilde F$, the boundedness of $B$, the convergence of Algorithm~\ref{alg:pdhgm} and the majorization above. The second is the rule of \cite[eq.~(19)]{ochs2015iteratively}, which checks the energy every $10$ iterations and stops at its first decrease. The run then stops by the rule of \cite[eq.~(20)]{ochs2015iteratively}, once $(J_\varphi(u_k)-J_\varphi(u_{k+1}))/J_\varphi(u_0)<10^{-7}$. Algorithm~\ref{alg:dca} (Ochs rules) in Section~\ref{sec:operators} uses the same two rules, with $J_\gamma$ and the subproblem \eqref{eq:subproblem}.

\paragraph{The method of Hintermüller and Wu.} Their Huberized model with $\mu=0$ uses the gradient $N\mathbb{K}$ on an $N\times N$ image and the Huber function $\phi_\gamma(s)=s^q/q-(1/q-1/2)\gamma^q$ for $s>\gamma$ and $\gamma^{q-2}s^2/2$ otherwise, with $\gamma=0.1$, the default of their experiments \cite[Section~4]{hintermuller2013nonconvex}, which we kept. With $\beta=\alpha N^q/q$ and the knot $g=\gamma/N$, it reads
\begin{equation*}
  J_{\rm HW}(u)=\tfrac12\|Au-f\|^2+\beta\sum_{i=1}^m\psi\bigl(\phi_i(u)\bigr),\qquad
  \psi(t)=\begin{cases}t^q-(1-\tfrac q2)g^q,& t>g,\\ \tfrac q2\,g^{q-2}t^2,& t\le g,\end{cases}
\end{equation*}
in our notation, and its value coincides with theirs. The function $\psi$ is convex on $[0,g]$, so it is not an admissible penalty, and their method has no box constraint.

The Newton matrices of \cite[(3.10), (3.14), (3.15), (3.20)]{hintermuller2013nonconvex} are the same in both sets of units, so we kept their parameters:
\begin{equation*}
\begin{gathered}
  \beta_{\max}=1.2-q,\quad c=1,\quad \rho_1=0.25,\quad \rho_2=0.75,\quad \kappa_1=0.25,\\
  \kappa_2=2,\quad \epsilon_d=10^{-8},\quad \tau_1=0.1,\quad \tau_2=0.9,
\end{gathered}
\end{equation*}
and $\varepsilon=10^{-4}\alpha=10^{-4}\beta qN^{-q}$~\cite[Algorithm~3.6]{hintermuller2013nonconvex}. As in their experiments, one inner iteration updates the regularization weight, and the Newton system is solved by conjugate gradients with residual tolerance $0.05$.

Denoising is the exception: they use a sparse Cholesky factorization there, whose cost cannot be counted in applications of $(A,A^\top)$. We used conjugate gradients with the same tolerance instead. This was the cheapest of the choices we tried: a tolerance of $10^{-10}$ lowers the objective at the stop by at most $0.05\%$ and takes five to seven times the work.

If conjugate gradients meet a direction of nonpositive curvature, we set the regularization weight to $\beta_{\max}$, as they do when they detect an indefinite matrix. The dual update \cite[(3.9)]{hintermuller2013nonconvex} uses the symmetrized blocks of their (3.13) with the projected dual. The Wolfe--Powell line search brackets an acceptable step and then bisects.

Their paper leaves two points open. It does not report the initial values, for which we took the regularization weight $\beta_{\max}$, the trust-region radius $1$ and the dual variable $p_0=\max(\phi(u_0),g)^{q-2}\mathbb{K}u_0$. It also does not say which regularization weight enters the dual update, and we used the one with which the Newton system was solved, not the one updated by the trust-region step.

\paragraph{Validation.} Table~\ref{tab:hw-validation} compares our implementation with \cite[Tables~5 and~6]{hintermuller2013nonconvex}. Their noise realization is not known, so we report means over three realizations for denoising and over five for tomography.

For denoising, the PSNR agrees to $0.1$~dB, the objective to $0.25\%$ and the sparsity to $0.006$ for every $q$. Two choices behind these values are ours. The phantom of our experiments does not have the gradient sparsity that they report for their test image, so we used a phantom of $512\times512$ pixels averaged over blocks of $2\times2$, which is close to it, as a guess at their image. For $q=0.25$, we used $\alpha=1.0\cdot10^{-3}$ in place of the value $1.2\cdot10^{-3}$ of their Figure~6(d), which did not reproduce their results.

For tomography, our PSNR is lower than theirs by $0.7$~dB for total variation and by $3.0$~dB for $q=0.75$, with their stopping tolerance $10^{-4}$. Our projector differs from theirs, the number of angles may differ as well (Appendix~\ref{app:repro}), and the PSNR depends strongly on the noise realization.

\begin{table}[tbp]
  \centering
  \caption{Our implementation of the method of \cite{hintermuller2013nonconvex} against the values of its Tables~5 (denoising) and~6 (tomography): the PSNR, the objective $J_{\rm HW}$ at the computed solution in their units, and the fraction of pixels with $\phi_i(u)>g$ (sparsity). Denoising on the block-averaged phantom, with $\alpha=1.0\cdot10^{-3}$ for $q=0.25$. Our values are means over the noise realizations.}
  \label{tab:hw-validation}
  \small
  \begin{tabular}{llrrrrrr}
    \toprule
    & & \multicolumn{2}{c}{PSNR} & \multicolumn{2}{c}{$J_{\rm HW}(\hat u_q)$} &
    \multicolumn{2}{c}{sparsity} \\
    \cmidrule(lr){3-4}\cmidrule(lr){5-6}\cmidrule(lr){7-8}
    problem & $q$ & ours & theirs & ours & theirs & ours & theirs \\
    \midrule
    denoising & TV & 37.47 & 37.57 & \textemdash & \textemdash & 0.192 & 0.196 \\
     & 0.75 & 40.91 & 41.00 & 134.02 & 134.30 & 0.057 & 0.058 \\
     & 0.5 & 40.97 & 41.02 & 116.46 & 116.67 & 0.053 & 0.053 \\
     & 0.25 & 39.96 & 39.93 & 113.85 & 114.00 & 0.056 & 0.050 \\
    \midrule
    tomography & TV & 26.09 & 26.80 & & & & \\
     & 0.75 & 34.43 & 37.39 & & & & \\
    \bottomrule
  \end{tabular}
\end{table}

\section{Reproducibility}
\label{app:repro}
\paragraph{Operators and data.} For deblurring, $A$ is the correlation with a Gaussian kernel normalized to unit sum. For tomography, $A$ is a parallel-beam X-ray transform that distributes each pixel over the two nearest detector bins with linear weights, with enough bins to cover the image diagonal and $15$ equispaced angles in $[0,\pi)$. This is our reading of the problem of \cite[Table~6]{hintermuller2013nonconvex}, which lists the angles $0,12,\dots,180$ degrees but reports a sinogram of $95\times13$. The transform and its transpose are computed by compiled loops over the same weights, so that $A^\top$ is the transpose of $A$ up to rounding errors, and $A$ is divided by an upper bound of its norm given by the Schur test, so that every operator satisfies $\|A\|\le1$. The noise is not clipped, and for tomography it is added before $A$ is divided by its norm bound, i.e., in the units of \cite{hintermuller2013nonconvex}; all methods receive the same $f$. The total variation starting point is computed from $\proj_B(f)$, or $\proj_B(A^\top f)$ for tomography, until its duality gap, relative to the energy at the start, falls below the tolerance of the corresponding nonconvex run. The step sizes $\tau=c_k\bar\tau$, $\sigma=\bar\tau/c_k$ of Algorithm~\ref{alg:dca} were used for all problems and for both IRL1 variants.

\paragraph{Test images.} The Kodak image kodim23 \cite{kodak} was converted to grayscale, cropped to its central $512\times512$ square and averaged over blocks of $2\times2$ pixels, so that no resampling is involved. For butterfly \cite{vanzalinge2019butterfly}, a grayscale square crop of $1184\times1184$ pixels was resized to $256\times256$ with an anti-aliasing filter.

\paragraph{Software.} All methods are implemented in Python and share the operators and
the data. The solvers depend on NumPy and on the package fastprox~\cite{fastprox}
(version 1.0.5), whose compiled routines compute the projections onto $B$ and onto the
dual balls \eqref{eq:dualball}. If Numba is installed, it compiles the loops of the blur
and of the X-ray transform, with the same results up to rounding errors. Scikit-image is
used to load and resize the test images, and matplotlib for the figures.

\end{appendices}

\backmatter

\section*{Declarations}
\bmhead{Competing interests} The authors have no relevant financial or non-financial interests to disclose.
\bmhead{Data availability} The test images are publicly available from the sources cited in Table~\ref{tab:setup}.
\bmhead{Code availability} The code that reproduces all experiments is archived on Zenodo~\cite{tvq_code}.
\bmhead{Use of AI tools} The conception of the method, including the DC splitting, its analysis, the algorithm and the design of the experiments, is the work of the authors. The authors implemented the core of the algorithm; starting from this implementation, Claude Opus 5.5 and Claude Fable 5.2 (Anthropic) were used to extend the code of the methods and the experiments, and to improve the presentation of the proofs and the text (Appendix~\ref{app:repro}). The authors verified the code and the proofs and take full responsibility for the content of this publication.

\bibliography{bibliografia,abstract_refs}

\end{document}